\documentclass[11pt]{amsart}

\usepackage{epigamath}

\usepackage[english]{babel}

\numberwithin{equation}{section}

\usepackage{amssymb}
\usepackage{verbatim}
\usepackage{amsmath,amsthm,mathrsfs}
\usepackage{amscd}
\usepackage{array}
\usepackage{tikz}
\usepackage[arrow,matrix]{xy}

\usepackage{bm}

\usepackage{stmaryrd}

\newtheorem{theo}{Theorem}[section]
\newtheorem{prop}[theo]{Proposition}
\newtheorem{lem}[theo]{Lemma}
\newtheorem{cor}[theo]{Corollary}
\newtheorem{theon}{Theorem}

\theoremstyle{definition}
\newtheorem{defi}[theo]{Definition}

\theoremstyle{remark}
\newtheorem{rem}[theo]{Remark}
\newtheorem{ex}[theo]{Example}
\newtheorem{con}[theo]{Convention}

\newcommand{\lser}[1]{(\!(#1)\!)}
\newcommand{\pow}[1]{\llbracket #1 \rrbracket}

\SetSymbolFont{stmry}{bold}{U}{stmry}{m}{n}

\def \si {{\sigma}}
\def \Ga {{\Gamma}}

\DeclareMathOperator{\Pic}{Pic}

\DeclareMathOperator{\Gal}{Gal}

\DeclareMathOperator{\Hom}{Hom}

\DeclareMathOperator{\GL}{GL}

\DeclareMathOperator{\nr}{nr}

\DeclareMathOperator{\sing}{sing}
\DeclareMathOperator{\red}{red}

\DeclareMathOperator{\Spec}{Spec}
\DeclareMathOperator{\cris}{cris}
\newcommand{\Char}{\operatorname{char}}
\DeclareMathOperator{\Bl}{Bl}
\DeclareMathOperator{\Proj}{Proj}
\DeclareMathOperator{\Spf}{Spf}

\def \P{{\mathbb P}}

\def\ov{\overline}

\def \Z {{\mathbb Z}}
\def \Q {{\mathbb Q}}
\def \F {{\mathbb F}}

\def \et {{\mathrm{\acute{e}t}}}

\DeclareMathOperator{\val}{val}

\DeclareMathOperator{\Kum}{Kum}

\def\G{{\mathbb G}}

\def\T{{\mathcal T}}

\def\sZ{{\mathcal Z}}
\def\sY{{\mathcal Y}}

\def\lra{\longrightarrow}

\def\H{{\rm H}}
\def\sA{{\mathcal A}}

\def\X{{\mathcal X}}

\def\O{{\mathcal O}}

\def\sY{{\mathcal Y}}

\def\si{\sigma}
\def\Ga{\Gamma}

\def\m{{\mathfrak m}}

\newcommand{\bthe}{\begin{theo}}
\newcommand{\ble}{\begin{lem}}
\newcommand{\bpr}{\begin{prop}}
\newcommand{\bco}{\begin{cor}}
\newcommand{\bde}{\begin{defi}}
\newcommand{\ethe}{\end{theo}}
\newcommand{\ele}{\end{lem}}
\newcommand{\epr}{\end{prop}}
\newcommand{\eco}{\end{cor}}
\newcommand{\ede}{\end{defi}}
\newcommand{\brem}{\begin{rem}}
\newcommand{\erem}{\end{rem}}
\newcommand{\bex}{\begin{ex}}
\newcommand{\eex}{\end{ex}}
\newcommand{\bcon}{\begin{con}}
\newcommand{\econ}{\end{con}}

\usepackage[shortlabels]{enumitem}
\setlist[enumerate,1]{label={\rm(\arabic*)}, ref={\rm\arabic*}} 

\EpigaVolumeYear{7}{2023} \EpigaArticleNr{10} \ReceivedOn{June 7, 2022}
\InFinalFormOn{November 9, 2022}
\AcceptedOn{December 14, 2022}

\title{Reduction of Kummer surfaces modulo 2 in the non-supersingular case}
\titlemark{Reduction of Kummer surfaces modulo 2 in the non-supersingular case}

\author{C.\,D.~Lazda}
\address{Department of Mathematics, Harrison Building, University of Exeter, EX4 4QF, United Kingdom}
\email{c.d.lazda@exeter.ac.uk}
\author{A.\,N.~Skorobogatov}
\address{Department of Mathematics, South Kensington Campus, Imperial College London, SW7 2AZ, United Kingdom\\
Institute for the Information Transmission Problems,
Russian Academy of Sciences,
Moscow, 127994, Russia}
\email{a.skorobogatov@imperial.ac.uk}

\authormark{C.\,D.~Lazda and A.\,N.~Skorobogatov}

\AbstractInEnglish{We obtain necessary and sufficient conditions for the good reduction of Kummer surfaces attached to abelian surfaces with non-supersingular reduction when the residue field is perfect of characteristic $2$. In this case, good reduction with an algebraic space model is equivalent to good reduction with a scheme model, which we explicitly construct.  We extend these results to Kummer surfaces attached to $2$-coverings of such abelian surfaces.}



\acknowledgement{The work on this paper started while A.S. was a fellow of FRIAS, University of Freiburg.  He is grateful to FRIAS for the hospitality and excellent working conditions.}

\begin{document}


\maketitle

\begin{prelims}

\DisplayAbstractInEnglish

\bigskip

\DisplayKeyWords

\medskip

\DisplayMSCclass

\end{prelims}


\newpage

\setcounter{tocdepth}{1}

\tableofcontents


\section{Introduction}

To an abelian surface $A$ over a field $K$, one can attach the associated Kummer surface $\Kum(A)$, defined as the minimal desingularisation of the quotient of $A$ by the antipodal involution $\iota(P)=-P$. If the characteristic of $K$ is not $2$, then $\Kum(A)$ is a K3 surface.  Quartic surfaces with sixteen nodes in $\P^3$ are classically studied singular models of Kummer surfaces associated to Jacobians of genus $2$ curves. When $K$ is a number field, such Kummer surfaces, as well as those attached to products of elliptic curves, form classes of K3 surfaces most accessible for studying their arithmetic; see for example \cite{HS16, SZ17} and references therein.

If the characteristic of $K$ is $2$, then the geometry of $\Kum(A)$ depends on the $2$-rank of $A$, that is, the integer $r\in \{0,1,2\}$ such that $A[2]$ has $2^r$ points over an algebraic closure of $K$. The abelian surface $A$ is called {\em ordinary} if $r=2$ and {\em supersingular} if $r=0$; in the remaining case, where $r=1$, $A$ will be called {\em almost ordinary}. Then $\Kum(A)$ is a K3 surface if and only if $A$ is not supersingular.  This was proved by Shioda when $A$ is a product of elliptic curves \cite{Shi74} and by Katsura in general \cite{Kat78}. In contrast, if $A$ is supersingular, then $\Kum(A)$ is geometrically rational. Katsura also points out another striking difference between these cases: if $A$ is not supersingular, then the singularities of $A/\iota$ are rational double points, which are four geometric points of type ${\rm D}_4$ when $A$ is ordinary, and two geometric points of type ${\rm D}_8$ when $A$ is almost ordinary, giving rise to sixteen rational curves in
$\Kum(A)$.  On the other hand, if $A$ is supersingular, then the unique singular point of $A/\iota$ is an elliptic singularity. In the non-supersingular case, the $2$-rank of $A$ also controls the height of the K3 surface $\Kum(A)$, which is $1$ in the ordinary case and $2$ in the almost ordinary case; see Remark~\ref{height}.

In this paper we are interested in necessary and sufficient conditions for good reduction of Kummer surfaces. We will therefore take $K$ to be a complete discretely valued field of characteristic zero, with ring of integers $\O_K$ and perfect residue field $k$. Given an abelian surface $A$ over $K$, we say that $X=\Kum(A)$ has {\em good reduction} if there exists a scheme or algebraic space $\X$ smooth and proper over $\O_K$ with generic fibre $\X_K\cong X$. Note that in this situation, the special fibre $\X_k$ is necessarily a K3 surface since it must have trivial canonical bundle $\omega_{\X_k}\cong \O_{\X_k}$,\footnote{This is because the isomorphism $\Pic(\X)\to \Pic(X)$ sends $\omega_{\X/\O_K}$ to $\omega_X$.}  and (coherent) Euler characteristic $\chi(\O_{\X_k})=2$. Although for general K3 surfaces, good reduction requires models which are algebraic spaces (see for example \cite[Example~5.2]{Mat15a}), in the cases considered in this paper, it turns out that schemes
suffice.

If $\Char(k)\neq 2$, then $\Kum(A)$ has good reduction if and only if there exists a quadratic twist $A^\chi$ of $A$ such that $A^\chi$ has good reduction. Indeed, if $\Kum(A)$ has good reduction, one can show that the inertia group of $K$ acts on ${\rm H}^1_{\et}(A_{\overline{K}},\Q_\ell)$ via a quadratic character and then apply the classical N\'eron--Ogg--Shafarevich criterion (see the proof of Theorem~\ref{theo: twist 1}, which is certainly well known to the experts). Conversely, any quadratic twist of $A$ satisfies $\Kum(A^\chi)\cong \Kum(A)$, so replacing $A$ with $A^{\chi}$, we can assume that $A$ has good reduction. In this case, the N\'eron model $\sA/\O_K$ of $A/K$ is an abelian scheme with generic fibre $\sA_K\cong A$. We can then form the quotient $\sA/\iota$ and blow up the singular subscheme to obtain a smooth model of $\Kum(A)$ over $\O_K$; see \cite[Lemma 4.2]{Mat15} and \cite[Proposition 3.11]{Ove21}.

If $\Char(k)=2$, one direction of this argument still works: if $\Kum(A)$ has good reduction, then a quadratic twist $A^\chi$ of $A$ has good reduction, and $\Kum(A)\cong \Kum(A^\chi)$. But the argument in the other direction breaks down, essentially because the singular subscheme of $\sA/\iota$ is no longer \'etale over $\O_K$. We are therefore left with the following question: suppose that $\Char(k)=2$, and let $A/K$ be an abelian surface with good reduction $\sA/\O_K$. When does $\Kum(A)$ have good reduction?

Our first result is a construction, in the non-supersingular case, of an explicit smooth model of $\Kum(A)$ over the ring of integers of a finite extension of $K$ that trivialises the Galois action on the $2$-torsion subgroup $A[2]$.  We use a local equation of the special fibre at each singular point, but not much else. Our method relies on the crucial fact that blowing up a (singular!) section commutes with specialisation to each fibre, provided the section meets both fibres at rational double points; see Proposition~\ref{lem: rdp blowups}.  As a consequence we obtain the following. 

\begin{theon} \label{theo: main 1} 
Assume that $\Char(k)=2$, and let $A/K$ be an abelian surface with good, non-supersingular reduction. Then $\Kum(A)$ has potentially good reduction with a scheme model.
\end{theon}

An important feature of our model is that it establishes a bijection between the sets of sixteen geometric exceptional curves of the special and generic fibres.\footnote{This bijection, however, does \emph{not} describe the images of these curves under the specialisation map on N\'eron--Severi groups.} A weaker result concerning good reduction with an algebraic space model, and without an explicit bound on the degree of the field extension required, follows from Artin's simultaneous resolution of singularities \cite{Art74}. This can be applied because $\sA/\iota$ is flat over~$\O_K$, and its fibres are normal varieties with at worst rational double points.

To explain more refined results than Theorem~\ref{theo: main 1}, let us fix an algebraic closure $\ov{K}$ of $K$, with residue field~$\bar{k}$, and let $\Ga_K$ denote the Galois group of $\ov{K}/K$. We therefore have the exact sequence
\begin{equation}
\label{eqn: exact sequence [2]} 0\lra \mathcal{A}[2]^\circ(\overline{K}) \lra A[2](\overline{K}) \lra \mathcal{A}[2](\bar{k})\lra 0
\end{equation}  
of $\Ga_K$-modules, where $\mathcal{A}[2]^\circ$ is the connected component of the identity of the $2$-torsion subscheme $\sA[2]\subset\sA$ (see Section~\ref{sec: prelim} below).

\begin{theon} \label{theo: main 2} Assume that $\Char(k)=2$, and let $A/K$ be an abelian surface with good, ordinary reduction. 
Then $\Kum(A)$ has good reduction over $K$ if and only if the exact sequence {\rm{(\ref{eqn: exact sequence [2]})}} of\, $\Ga_K$-modules splits. Moreover, in this case $\Kum(A)$ has good reduction with a scheme model.
\end{theon}

Another way of phrasing the given condition is as the splitting of the connected-\'etale sequence
\[ 0\lra \sA[2]^\circ \lra \sA[2] \lra \sA[2]^\et\lra 0 \]
of finite flat group schemes over $\O_K$. It is a consequence of Cartier duality that the splitting of {\rm{(\ref{eqn: exact sequence [2]})}} implies in particular that the $\Ga_K$-module $A[2](\overline{K})$ is unramified (see Section~\ref{sec: prelim}); however, the converse fails~-- there exist examples where $A[2](\overline{K})$ is unramified but the exact sequence {\rm{(\ref{eqn: exact sequence [2]})}} is non-split (see Example~\ref{exa: nonsplit}). In~particular, this gives an example of a K3 surface over $\Q_2$ which only attains good reduction over a non-trivial unramified extension.\footnote{The only examples of such behaviour in the literature are for K3 surfaces over $\Q_p$ with $p\geq 5$; see \cite[Theorem 7.2]{LM18}.}

We also have similar results in the almost ordinary case.

\begin{theon} \label{theo: main 3} Assume that $\Char(k)=2$, and let $A/K$ be an abelian surface with good, almost ordinary reduction. Then $\Kum(A)$ has good reduction over $K$ if and only if the $\Ga_K$-module $A[2](\overline{K})$ is trivial. Moreover, in this case $\Kum(A)$ has good reduction with a scheme model.
\end{theon}

Again, it is easy to see that this condition is not vacuous: there exist examples where the $\Ga_K$-module $A[2](\overline{K})$ is unramified but non-trivial.

The `if' directions of Theorems~\ref{theo: main 2} and~\ref{theo: main 3} are proved by explicitly constructing a smooth model for $\Kum(A)$, by resolving the singularities of the quotient scheme $\sA/\iota$ as in the proof of Theorem~\ref{theo: main 1}. The hypotheses on $A[2](\overline{K})$ are exactly what is required to make this construction work over $\O_K$ itself. 

The `only if' directions of Theorems~\ref{theo: main 2} and~\ref{theo: main 3} are rather more involved, and the proofs use (the easy direction of) the main result of \cite{CLL}. (Note that since we have potentially good reduction by Theorem~\ref{theo: main 1}, the crucial Hypothesis~($\star$) of \cite{CLL} is always satisfied in our case.)  Indeed, we have an isomorphism $\sA_k/\iota\cong(\sA/\iota)_k$ (see Proposition~\ref{spec quot}), and this implies the key observation that $\Kum(\sA_k)$ is the `canonical reduction' of $\Kum(A)$ in the terminology introduced in \cite{CLL} (see Section~\ref{sec: grc} for the definition). By \cite[Theorem 1.6]{CLL}, good reduction of $\Kum(A)$ with an algebraic space model is therefore equivalent to the existence of an isomorphism between the Galois representations in the $\ell$-adic \'etale cohomology groups of $\Kum(\sA_k)$ and $\Kum(A)$, where $\ell$ is any odd prime.  We compare the Galois action on these cohomology groups by explicitly calculating
the Galois action on the exceptional curves of $\Kum(\sA_k)$ and $\Kum(A)$.

We finish this paper by applying our method to give a necessary and sufficient condition for good reduction of `twisted' Kummer surfaces, that is, Kummer surfaces associated to $2$-coverings of abelian surfaces; see Theorem~\ref{t2}.  Again, we only have results in the case of non-supersingular reduction. Indeed, the approach of this paper does not seem to be suitable for studying good reduction of Kummer surfaces attached to abelian surfaces with good, supersingular reduction. In the non-supersingular case, we make crucial use of the fact that a smooth model for $\Kum(A)$ over $\O_K$ can be obtained by resolving the singularities of the obvious singular model $\sA/\iota$.  This completely breaks down in the supersingular case, where no such resolution of $\sA/\iota$ can possibly provide a smooth model for $\Kum(A)$.

\subsection*{Notation and conventions}

Throughout, $\O_K$ is a complete discrete valuation ring with field of fractions $K$, maximal ideal $\m_K$, and residue field $k=\O_K/\m_K$. We assume that $\Char(K)=0$ and that $k$ is perfect with $\Char(k)=2$.

Let $\ov K$ be an algebraic closure of $K$, and write $\Ga_K=\Gal(\ov K/K)$.  The residue field of the maximal unramified extension $K_{\nr}\subset\ov K$ is an algebraic closure $\bar k$ of $k$. We write $\Ga_k=\Gal(\bar k/k)=\Gal(K_{\nr}/K)$, and we have the usual inertia exact sequence
\[ 1 \lra {\rm I}_K \lra \Gamma_K \lra \Ga_k \lra 0\]
relating $\Ga_K$ and $\Ga_k$.

Let $G$ be the (abstract) group $\Z/2$. For any abelian group object $A$ (in any category), we write $\iota\colon A\rightarrow A$ for the involution $\iota(x)=-x$, which we also think of as an action of $G$.

\subsection*{Acknowledgments}

This paper benefited from intensive discussions with Evis Ieronymou to whom we are very grateful.  We would like to thank Martin Liebeck for suggesting to use group cohomology in the proof of Theorem~\ref{t0}. We thank Otto Overkamp for his interest in this paper and helpful remarks.  

\section{Geometry in the ordinary case} \label{sec: geom ord}

For this section, we will let $A/k$ be an ordinary abelian variety. Thus $A[2](\bar k)\cong (\Z/2)^2$, and the quotient variety $A/G$ has four (geometric) singular points, all of type ${\rm D}_4^1$; see \cite[table on p.~144]{Sch09}. These points are in a natural $\Ga_k$-equivariant bijection with $A[2](\bar k)$. The exceptional divisor of the minimal resolution $\Kum(A)\to A/G$ therefore contains sixteen rational curves, occurring in four (geometrically) connected components. Our goal here is to calculate the Galois action on these curves and in particular prove the following result.

\begin{prop} \label{prop: galois exceptional ordinary} The sixteen exceptional curves of\, $\Kum(A)\rightarrow A/G$ are indexed\, $\Ga_k$-equivariantly by the set $A^\vee[2](\bar k)\times A[2](\bar k)$.
\end{prop}

We achieve this by first considering the case where $k$ is algebraically closed and giving a geometric description of how these singular points arise. This description will then be easily seen to be Galois equivariant. Note that the factor $A[2](\bar k)$ arises from the four singular points of $A/G$. It is therefore enough to complete around the image of the identity $O\in A$ and prove that the exceptional curves of the minimal resolution of the resulting ${\rm D}_4^1$ singularity are naturally indexed by $A^\vee[2](\bar k)$.

So suppose for now that $k$ is algebraically closed, and let $\widehat{\O}_{A,O}$ be the completed local ring of $A$ at the identity element $O\in A$.

\begin{prop} \label{prop: ord prod} There exists an isomorphism
\[ \widehat{\O}_{A,O}\cong k\pow{u,v}\]
such that the involution $\iota$ acts via $\iota(u)=\frac{u}{1+u}$ and $\iota(v)=\frac{v}{1+v}$, and the formal group law $\Delta$ satisfies
\[ \Delta(u)\equiv 1\otimes u + u\otimes 1 + u\otimes u \mod\;(u^2\otimes 1,1\otimes u^2,v^2\otimes 1, 1\otimes v^2 ),\]
\[ \Delta(v)\equiv 1\otimes v + v\otimes 1 + v\otimes v \mod\;(u^2\otimes 1,1\otimes u^2,v^2\otimes 1, 1\otimes v^2 ).\]
\end{prop} 

\begin{proof}
As in the proof of \cite[Lemma 4]{Kat78}, the classification of formal group laws over $k$ implies 
that there is an isomorphism of formal groups
\begin{equation}
\widehat{A}_{/O}\cong \widehat{E\times E}_{/O},\label{formal}
\end{equation}
where $E$ is any ordinary elliptic curve over $k$. For example, we can take $E$ to be the curve
$$y^2+xy=x^3+x.$$ Note that the involution $\iota$ on $E$ is defined by $\iota(x,y)=(x,y+x)$.  The function $z=x/y$ is a local parameter at the origin of the group law $O\in E$, so that the completed local ring $\widehat{\mathcal{O}}_{E,O}$ is isomorphic to $k\pow{z}$.  We have $\iota(z)=z(1+z)^{-1}$. In view of (\ref{formal}), this proves the claim concerning the action of $\iota$. For the claim concerning the form of $\Delta$, we simply use the formula on the top of p.~115 of \cite{Sil86}.
\end{proof}

Let $S=\Spec(\widehat{\mathcal{O}}_{A,O})$, let $T=S/G$, and let $q\colon S\to T$ be the quotient morphism.  It follows from \cite[Proposition 1.1]{Sch09}, which is simply a statement of the main result of \cite{Art75}, that
$$T=\Spec\left(\frac{k\pow{x,y,z}}{(z^2+xyz+x^2y+xy^2)}\right),$$
where $q$ is given by
\[ x=\frac{u^2}{1+u},\;\;\;\;y=\frac{v^2}{1+v},\;\;\;\; z=\frac{uv(u+v)}{(1+u)(1+v)}.\]
Let $O\in S$ denote the closed point and $q(O)\in T$ its image in $T$; thus $O=V(u,v)$ and $q(O)=V(x,y,z)$ as closed subschemes of $S$ and $T$, respectively. Let
\[ S^{(1)}:=\Bl_OS\lra S,\;\;\;\; T^{(1)}=\Bl_{q(O)}T\lra T\]
be the blowups at $O$ and $q(O)$, respectively. Explicitly, we have
\[ S^{(1)}=\Proj\left(\frac{k\pow{u,v}[U,V]}{(uV+vU)}\right)\]
and 
\[ T^{(1)}= \Proj\left(\frac{k\pow{x,y,z}[X,Y,Z]}{(xY+yX,xZ+zX,yZ+zY,Z^2+xYZ+xXY+xY^2)}\right). \]
It is clear that $S^{(1)}$ is smooth. By \cite[Proposition 3(i)]{Kat78} or \cite[Proposition 5.1]{Sch09}, the singular point $q(O)$ is a rational double point of type ${\rm D}_4^1$ (for the classification of rational double points in all characteristics, see \cite{Art77}). By \cite[Corollary 1.7]{Sch09}, $T^{(1)}$ is normal with at most rational double points as singularities, and the minimal desingularisation of $T$ factors through $T^{(1)}$. We let $D^{(1)}\subset S^{(1)}$ and $E^{(1)}\subset T^{(1)}$ denote the reduced exceptional subschemes of $S^{(1)}\to S$ and $T^{(1)}\to T$, respectively.

\begin{lem}  The rational map $q\colon S^{(1)}\dashrightarrow T^{(1)}$ is regular, and the induced map $D^{(1)}\rightarrow E^{(1)}$ is a universal homeomorphism.
\end{lem}

\begin{proof}
Let $S'\to S$ be the blowup of the fibre $q^{-1}(q(O))\subset S$.  By the universal property of blowing up, we have a commutative diagram
$$
\xymatrix{S'\ar[r]\ar[d]&S\ar[d]\\
T^{(1)} \ar[r]&T\rlap{.}}
$$
By \cite[Lemma 1.3, Remark 1.4]{Sch09}, the ideal of $q^{-1}(q(O))$ is the bracket ideal $\m_{O}^{[2]}=\{f^2\mid f\in\m_O\}$, where $\m_{O}\subset\widehat{\mathcal{O}}_{A,O}$ is the maximal ideal. We have $\m_{O}^{[2]}\,\m_{O}=\m_{O}^3$, which is simply saying that $(u^2,v^2)(u,v)=(u,v)^3$.  Thus the ideal $\m_{O}^{[2]}$ becomes invertible on $S^{(1)}$. By the universal property of blowing up, there is a unique factorisation $S^{(1)} \to S'\to S$, which therefore provides the required regular map $S^{(1)}\rightarrow T^{(1)}$ extending $q$. Explicitly, the map $q$ is given by
\[ (u,v,[U:V])\longmapsto \left(\frac{u^2}{1+u},\frac{v^2}{1+v},\frac{uv(u+v)}{(1+u)(1+v)},\left[ \frac{U^2}{1+u}:\frac{V^2}{1+v}:\frac{UV(u+v)}{(1+u)(1+v)}\right]\right), \]
which on $D^{(1)}=V(u,v)\subset S^{(1)}$ simplifies to 
\[ (0,0,[U:V])\longmapsto (0,0,0,[U^2:V^2:0]).\]
This map can therefore be identified with the relative Frobenius and is thus a universal homeomorphism.
\end{proof}

It follows from the results of \cite{Sch09} that $T^{(1)}$ has precisely three singular points, all of type ${\rm A}_1$. An explicit calculation gives these three points as 
\[ (x,y,z,[X:Y:Z])=(0,0,0,[0:1:0]),\;\;(0,0,0,[1:0:0]),\;\;(0,0,0,[1:1:0]). \]
Since $D^{(1)}\to E^{(1)}$ is a universal homeomorphism and $k$ is algebraically closed, these in turn give rise to three $k$-points on $D^{(1)}$. 
Explicitly, these are the three points
\[ (u,v,[U:V])=(0,0,[1:0]),\;\;(0,0,[0:1]),\;\;(0,0,[1:1]).\]
Via the natural identification $D^{(1)}=\P(T_OA)$, where $T_O A$ is the tangent space to $A$ at $O$, these three points therefore correspond to three distinguished tangent directions to $A$ at $O$. We can give an alternative description of these tangent directions as follows.

Since $k$ has characteristic $2$, the inclusion $A[2]\subset A$ induces an isomorphism $T_O(A[2])\,\tilde\to\, T_OA$, and the formula for $\Delta$ given in Proposition~\ref{prop: ord prod} gives
$A[2]^\circ \cong \mu_2^{\oplus 2}$. We therefore obtain three subgroups of $A[2]^\circ$ isomorphic to $\mu_2$, and the tangent spaces to these three subgroups give three preferred tangent directions to $A$ at $O$. 

\ble 
These three directions are precisely those coming from the three singular points of\, $T^{(1)}$.
\ele

\begin{proof}
Again, we can just calculate everything explicitly. 
Indeed, $A[2]^\circ$ is the closed subgroup scheme of $S$
defined by $u=\iota(u)$, $v=\iota(v)$; thus we can write 
\[ A[2]^\circ \cong \Spec\left(\frac{k\pow{u,v}}{(u^2,v^2)}\right) \]
with the group law given by $\Delta(u)=1\otimes u + u\otimes 1 + u\otimes u$, and similarly for $v$. The three preferred tangent directions come from the three closed subgroup schemes of $\Spec\left(\frac{k\pow{u,v}}{(u^2,v^2)}\right)$ defined by
\[ u=0,\;\;\;\;v=0,\;\;\;\;u+v=0, \]
respectively. The strict transform of the closed subscheme $\{u=0\}\subset S$ is the closed subscheme $\{U=0\}$ of $S^{(1)}$, which intersects the exceptional divisor $D^{(1)}$ at the point $(u,v,[U:V])=(0,0,[0:1])$. Similarly, the strict transform of the closed subscheme $\{v=0\}\subset S$ (respectively, $\{u+v=0\}$)  intersects the exceptional divisor $D^{(1)}$ at the point $(0,0,[1:0])$ (respectively, $(0,0,[1:1])$).
\end{proof}

We now drop the assumption that $k$ is algebraically closed. Then we can form $S=\Spec(\widehat{\mathcal{O}}_{A,O})$, $T:=S/G$, $q\colon S\rightarrow T$, and $T^{(1)}:=\Bl_{q(O)}T$ as above, and thus $T^{(1)}$ has precisely three singular points over $\bar{k}$. The description of these three singular points as arising from the three subgroups of $A_{\bar k}[2]^\circ$ respects the Galois action, and so we obtain the following.

\bco 
\label{2.1}
Let $A$ be an ordinary abelian surface over $k$. Then
there is an isomorphism of $k$-schemes $(T^{(1)}_{\sing})_{\red}\cong (A^\vee[2])_{\red}\setminus \{O\}$.
\eco

\begin{proof}
The three subgroups of $A_{\bar k}[2]^\circ$ isomorphic to $\mu_2$ canonically
correspond to the three non-zero elements of the abelian group
$\Hom(A^\vee[2](\bar k),\F_2)$, which we can identify canonically with $A^\vee[2](\bar k)$ via the isomorphism $\wedge^2 (A^\vee[2](\bar k))\cong \F_2$. These isomorphisms are $\Ga_k$-equivariant.
\end{proof} 

If we now consider the minimal resolution $\Kum(A)\rightarrow A/G$, then we can apply the description given in Corollary~\ref{2.1} around each of the four singular points of $A/G$. Since these singular points correspond precisely to the points of $A[2](\bar k)$, this completes the proof of Proposition~\ref{prop: galois exceptional ordinary}.\qed

\section{Geometry in the almost ordinary case} \label{sec: geom almost}

We now want to carry out a similar analysis when $A$ is an almost ordinary abelian surface over $k$. We have $A[2](k)= A[2](\bar{k})=\Z/2$, and  $A/G$ has two geometric singular points, both $k$-rational and both of type ${\rm D}_8^2$; see  \cite[table on p.~144]{Sch09}. The exceptional divisor of the minimal resolution $\Kum(A)\rightarrow A/G$ contains sixteen rational curves, occurring in two (geometrically) connected components. The main result of this section is then the following.

\begin{prop} \label{prop: rational curves} The sixteen exceptional curves of\, $\Kum(A)\rightarrow A/G$ are all $k$-rational.
\end{prop}

Translating by the unique non-identity point of $A[2](k)= A[2](\bar{k})=\Z/2$, we see that it suffices to prove the rationality of the eight exceptional curves of the minimal resolution of the formal ${\rm D}^2_8$ singularity $\Spec(\widehat{\mathcal{\O}}_{A,O})/G$. If we look at the associated Dynkin diagram, we see that the only possible Galois action is to interchange the two curves corresponding to the nodes at the two `short ends' of the diagram. Thus we see straight away that at least six of the exceptional curves are $k$-rational, and the remaining two are defined (at worst) over a quadratic extension of $k$.

\subsection{Formal groups and \texorpdfstring{$\bm{2}$}{2}-divisible groups} \label{sec: formal split}

To prove Proposition~\ref{prop: rational curves}, we decompose the formal group of $A$.

\begin{lem} There is a $k$-isomorphism of formal groups $\Spf(\widehat{\mathcal{O}}_{A,O})\cong G_1\times_k G_2$, where $G_h$ is a connected, $2$-divisible formal group of dimension $1$ and height $h$, for $h=1,2$.
\end{lem}

\begin{proof}
Let $A[2^\infty]$ be the $2$-divisible group of $A$ and $A^\vee[2^\infty]$ that of its dual. If we let $\mathbf{D}$ denote the Cartier duality functor for $p$-divisible groups (so that $\mathbf{D}(\mathcal{G})[p^n]=\Hom(\mathcal{G}[p^n],\mu_{p^n})$), then the Weil pairing gives an isomorphism $A^\vee[2^\infty] \cong \mathbf{D}(A[2^\infty])$. Since $k$ is perfect, the connected-\'etale sequences of both $A[2^\infty]$ and $A^\vee[2^\infty]$ split canonically by \cite[last paragraph on p.~142]{Tat97}. Thus by Cartier duality, we obtain a decomposition
\[ A[2^\infty] \cong \mathbf{D}\left(A^\vee[2^\infty]^{\et}\right) \times \mathcal{G}_2 \times A[2^\infty]^{\et}, \]
where $\mathcal{G}_2$ is a connected $2$-divisible group of height $2$ and dimension $1$. Passing to the connected component of $A[2^\infty]$, we find that
\[ A[2^\infty]^\circ \cong \mathbf{D}\left(A^\vee[2^\infty]^{\et}\right) \times \mathcal{G}_2 ,\]
and invoking \cite[Proposition 1]{Tat67} then gives the result. 
\end{proof}

\begin{lem}
Let $G_h=\Spf(k\pow{u})$ be a $2$-divisible formal group over $k$ of dimension $1$ and height $h\geq 1$. Then the inversion map $\iota$ is of the form
\[ \iota(u) = u+a \]
for some $a\in k\pow{u}^\times$ with $\val_u(a)=2^h$.
\end{lem}

\begin{proof}
If $\iota(u)=a_1u+a_2u^2+\cdots $, then we must have $a_1=1$ and hence $\iota(u)=u+a$ for some $a\in k\pow{u}$ with $\val_u(a)\geq 2$. The $2$-torsion $G_h[2]$ is the $\iota$-fixed subscheme of $G_h$, which is defined by the ideal generated by $u+\iota(u)=a$ and is therefore of rank $\val_u(a)$ over $k$. But since $G_h$ has height $h$, we know that $G_h[2]$ is of rank $2^h$ over $k$, and the result follows.   
\end{proof}

\begin{cor} \label{cor1}
Let $x=u\cdot\iota(u)$ and $a=u+\iota(u)$. 
The quotient $G_h/\iota$ is given by $\Spf(k\pow{x})$, the element $a$ lies in $k\pow{x}$ and satisfies $\val_x(a)=2^{h-1}$, and $u$ satisfies the equation
\[ u^2+au+x=0. \]
\end{cor}

\begin{proof}
Arguing as in the proof of \cite[Lemma 1]{Art75}, we see that $k\pow{u}$ is finite free over $k\pow{x}$ of rank $\dim_k k\pow{u}/(x)=2$. In particular, this implies that $k\lser{u}/k\lser{x}$ is finite of degree $2$, from which we deduce that the inclusion $k\lser{x}\subset k\lser{u}^G$ is an equality, and hence $k\pow{u}^G=k\lser{x}\cap k\pow{u}=k\pow{x}$. Then $a=u+\iota(u)$ is an element of $k\pow{u}^G=k\pow{x}$. Moreover, $\val_u(a)=\val_u(u^{2^h}f)=2^h$ and $\val_u(x)=2$; hence $\val_x(a)=2^{h-1}$, as required. Finally the fact that $u^2+au+x=0$ is a straightforward check.
\end{proof}

\subsection{A na\"{i}ve resolution of \texorpdfstring{\bm{$A/G$}}{A/G}} \label{sec: naive}

Having decomposed the formal group of $A$, we can then explicitly resolve the singular surface $\Spec(\widehat{\mathcal{\O}}_{A,O})/G$ in the most na\"ive way, by repeatedly blowing up the `worst' singularity. Let us write $S=\Spec(\widehat{\O}_{A,O})=\Spec(k\pow{u,v})$, where we may use the results of Section~\ref{sec: formal split} above to choose $u$ and $v$ in such a way that $G$ acts via
\[ \iota(u)=u+u^2f(u),\;\;\;\;\iota(v)=v+v^4g(v) \]
for some $f\in k\pow{u}^\times$, $g\in k\pow{v}^\times$. 
Let $x=u\iota(u)$, $y=v\iota(v)$, $z=u\iota(v)+v\iota(u)$.
From Corollary~\ref{cor1}, we obtain $u^2+au+x=0$, where $a=x\epsilon(x)$ for some
$\epsilon(x)\in k\pow{x}^\times$. Similarly, we have $v^2+bv+y=0$, where
$b=y^2\eta(y)$ for some $\eta(y)\in k\pow{y}^\times$.
By \cite[Proposition 1.1]{Sch09}, we see that the quotient $T:=S/G$ is given by 
\[ T=\Spec\left(\frac{k\pow{x,y,z}}{(z^2+xy^2z\epsilon(x)\eta(y)+x^2y\epsilon(x)^2+xy^4\eta(y)^2)}\right).\]
Write $q\colon S\to T$ for the quotient morphism.
Then $q(O)$ is the unique singular point of $T$ and is of type ${\rm D}_8^2$.

\begin{lem} There exist a cubic extension $k'/k$ and a change of co-ordinates $x\mapsto \alpha x$, $y\mapsto \beta y$, $z\mapsto \gamma z$, with $\alpha,\beta,\gamma\in k'\pow{x,y,z}^\times$, inducing an isomorphism
\[ T\cong \Spec\left(\frac{k\pow{x,y,z}}{(z^2+xy^2z+x^2y+xy^4)}\right).\]
\end{lem}

\begin{proof}
If we make the given substitution, then we obtain the equation
\[ \gamma^2z^2+\alpha\beta^2\gamma xy^2z\epsilon\eta+\alpha^2\beta x^2y\epsilon^2+\alpha\beta^4xy^4\eta^2=0, \]
and for this to be of the required form, we need
\[ \alpha\beta^2\epsilon\eta=\gamma,\;\;\;\;\alpha^2\beta\epsilon^2=\gamma^2,\;\;\;\;\alpha\beta^4\eta^2=\gamma^2. \]
Solving these equations gives
\[ \alpha=\epsilon^{-2},\;\;\;\;\beta^3=\eta^{-2},\;\;\;\;\gamma=\alpha\beta^2\epsilon\eta. \]
Thus we need to be able to take a cube root of $\eta\in k\pow{y}^\times$, which we can always do over the (at worst) cubic extension $k(\sqrt[3]{\eta(0)})$.
\end{proof}

We already know that all the exceptional curves in the minimal resolution of $T$ are rational over a quadratic extension of $k$. To prove that they are $k$-rational, we are therefore allowed to make a cubic extension of $k$ and hence assume that
\[ T\cong \Spec\left(\frac{k\pow{x,y,z}}{(z^2+xy^2z+x^2y+xy^4)}\right).\]
Let $T^{(1)}$ be the blowup of $T$ at $q(O)=V(x,y,z)$. Thus
\[ T^{(1)}= \Proj\left(\frac{k\pow{x,y,z}[X,Y,Z]}{(xY+yX,xZ+zX,yZ+zY,Z^2+y^2XZ+yX^2+y^3XY)} \right). \]
Let $E^{(1)}\subset T^{(1)}$ denote the reduced exceptional subscheme of $T^{(1)}\rightarrow T$.

Via direct calculation, we can show that $T^{(1)}$ has precisely two singular points, both defined over $k$, and given in co-ordinates $(x,y,z,[X:Y:Z])$ by 
\[ Q^{(1)}_1=(0,0,0,[1:0:0]),\;\;\;\;Q_2^{(1)}=(0,0,0,[0:1:0]). \]
These are of type ${\rm A}_1$ and ${\rm D}_6^1$, respectively. Indeed, if we set $y'=y/x$ and $z'=z/x$, we can explicitly describe the formal completion $T^{(1)}_{/Q_1^{(1)}}$ as
\[T^{(1)}_{/Q_1^{(1)}} = \Spec\left( \frac{k\pow{x,y',z'}}{(z'^2 + x^2y'^2z' + xy'  + x^3y'^4)}\right), \]
and if we set $x'=x/y$ and $z''=z/y$, then the formal completion $T^{(1)}_{/Q_2^{(1)}}$ is given by
\[T^{(1)}_{/Q_2^{(1)}} = \Spec\left( \frac{k\pow{x',y,z''}}{(z''^2+x'y^2z''+x'^2y+x'y^3)}\right). \]
We now let $T^{(2)}\rightarrow T^{(1)}$ be the blowup at the $k$-rational singular point $Q^{(1)}_2$ of type ${\rm D}_6^1$. We let $E^{(2)}\subset T^{(2)}$ denote the reduced exceptional subscheme of $T^{(2)}\rightarrow T$, which naturally decomposes as $E^{(2)}=E^{(2)}_1\cup E^{(2)}_2$, where $E^{(2)}_1$ is the (reduced) strict transform of $E^{(1)}$ and $E^{(2)}_2$ is the reduced exceptional subscheme of $T^{(2)}\rightarrow T^{(1)}$. 

To perform explicit calculations on $T^{(2)}$, we write $x'=x/y$ and $z''=z/y$ and take the formal completion
\[ T^{(1)}_{/Q_2^{(1)}}=\Spec\left(\frac{k\pow{x',y,z''}}{(z''^2+x'y'^2z''+x'^2y+x'y^3)}\right) \]
at $Q_2^{(1)}$. Thus $T^{(2)}\times_{T^{(1)}} T^{(1)}_{/Q_2^{(1)}}$ is given by
\[\Proj\left(\frac{k\pow{x',y,z''}[X',Y',Z']}{\left(x'Y'+yX',x'Z'+z''X',yZ'+z''Y',
Z'^2+y'^2X'Z'+yX'^2+yX'Y'\right)} \right). \]
Further explicit calculations show that $T^{(2)}$ has three singular points $Q_1^{(2)}$, $Q_2^{(2)}$, and $Q_3^{(2)}$. The first of these, $Q_1^{(2)}$, is obtained simply as the strict transform of $Q_1^{(1)}$, lies on $E_1^{(2)}\setminus E_2^{(2)}$, and is $k$-rational of type ${\rm A}_1$. The second, $Q_2^{(2)}$, is the unique point of $E_1^{(2)}\cap E_2^{(2)}$ and is given in co-ordinates $(x',y,z'',[X':Y':Z'])$ by $(0,0,0,[1:0:0])$. It is therefore $k$-rational, and it is not difficult to show, again via explicit calculation, that it is of type ${\rm A}_1$.

The final singularity $Q^{(2)}_3$ lies on $E_2^{(2)}\setminus E_1^{(2)}$. In co-ordinates $(x',y,z'',[X':Y':Z'])$, it is given by $(0,0,0,[0:1:0])$, and if we set $x''=x'/y$ and $z'''=z''/y$, then the formal completion $T^{(2)}_{/Q^{(2)}_3}$ is given explicitly by 
\[ \Spec\left( \frac{k\pow{x'',y,z'''}}{(z'''^2+x''y^2z'''+x''^2y+x''y^2)} \right).\]
In particular, this singularity is $k$-rational and of type ${\rm D}_4^0$.

We now let $T^{(3)}\rightarrow T^{(2)}$ be the blowup at $Q^{(2)}_3$. We let $E^{(3)}_3\subset T^{(3)}$ denote the reduced exceptional subscheme of $T^{(3)}\rightarrow T^{(2)}$, write $E^{(3)}_1$, $E^{(3)}_2$ for the (reduced) strict transforms of $E^{(2)}_1$, $E^{(2)}_2$, and set $E^{(3)}=E^{(3)}_1\cup E^{(3)}_2\cup E^{(3)}_3$. A further explicit calculation shows that $T^{(3)}$ has five singular points, all of type ${\rm A}_1$, and lying on $E^{(3)}$ as follows: 
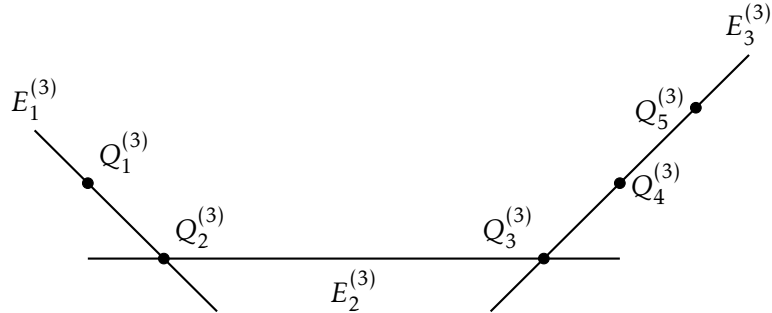
\begin{center}
\begin{tikzpicture}
\draw[thick] (-1,0) -- (6,0);
\draw[thick] (.7,-.7) -- (-1.7,1.7);
\draw[thick] (4.3,-.7) -- (7.7,2.7);
\filldraw[black] (0,0) circle (2pt) node[anchor=south west]{$Q_2^{(3)}$};
\filldraw[black] (5,0) circle (2pt) node[anchor=south east]{$Q_3^{(3)}$};
\filldraw[black] (6,1) circle (2pt) node[anchor=west]{$Q_4^{(3)}$};
\filldraw[black] (7,2) circle (2pt) node[anchor=east]{$Q_5^{(3)}$};
\filldraw[black] (-1,1) circle (2pt) node[anchor=south west]{$Q_1^{(3)}$};
\filldraw[black] (-1.7,1.7) circle (0pt) node[anchor=south]{$E_1^{(3)}$};
\filldraw[black] (7.7,2.7) circle (0pt) node[anchor=south]{$E_3^{(3)}$};
\filldraw[black] (2.5,0) circle (0pt) node[anchor=north]{$E_2^{(3)}$};
\end{tikzpicture}
\end{center}
It is clear that the points $Q_1^{(3)}$ and $Q_2^{(3)}$ are $k$-rational since they are simply the unique preimages of $Q_1^{(2)}$ and $Q_2^{(2)}$. That $Q_3^{(3)}$ is $k$-rational then follows from the geometry of the above configuration since there is no other point it could be sent to under the Galois action. To see that the last two points $Q_4^{(3)}$ and $Q_5^{(3)}$ are $k$-rational, we can do an explicit  calculation. Indeed, if we write
\[ T^{(2)}_{/Q^{(2)}_3}=\Spec\left( \frac{k\pow{x'',y,z'''}}{(z'''^2+x''y^2z'''+x''^2y+x''y^2)} \right),
\]
so that $T^{(3)}\times_{T^{(2)}} T^{(2)}_{/Q^{(2)}_3}$ is given by
\[ \Proj\left(\frac{k\pow{x'',y,z'''}[X'':Y'':Z'']}{(x''Y''+yX'',x''Z''+z'''X'',yZ''+z'''Y'',Z''^2+y^2X''Z''+yX''^2+yX''Y'')}\right), \]
then we can explicitly give the singular points $Q_3^{(3)}$, $Q_4^{(3)}$, and $Q_5^{(3)}$ in co-ordinates as
\begin{align*} (x'',y,z''',[X'':Y'':Z'']) =(0,0,0,[1:0:0]),\;\; 
(0,0,0,[1:1:0]), \;\;
 (0,0,0,[0:1:0]).
 \end{align*}
This completes the proof of Proposition~\ref{prop: rational curves}.\qed

\section{The Galois action on Kummer surfaces} 

We can now use the above results to describe, in the non-supersingular case, the Galois action on the cohomology of Kummer surfaces over $K$ and $k$. We will assume that $k$ is perfect of characteristic $2$ and let $\sA/\O_K$ be a relative abelian surface with generic fibre $A:=\sA_K$.

 \subsection{Preliminaries on Cartier duality} \label{sec: prelim}

The $2$-torsion subscheme $\sA[2]$ is a finite flat group scheme over $\mathcal{O}_K$; we let $\sA[2]^\circ$ be the connected component
of the identity. By  
\cite[Section~(3.7)(I)]{Tat97},
we have the connected-\'etale
exact sequence of finite flat group $\mathcal{O}_K$-schemes
\begin{equation*}
0\lra \sA[2]^\circ\lra \sA[2]\lra \sA[2]^\et\lra 0,
\end{equation*}
where $\sA[2]^\et$ is the largest \'etale group $\mathcal{O}_K$-scheme quotient of $\sA[2]$. If we base change from $\mathcal{O}_K$ to $k$, the connected-\'etale sequence splits, and we have
\[ \left(\sA[2]^\et\right)_k=(\sA_k[2])^\et=(\sA_k[2])_{\red},\]
and this is the unique \'etale group scheme over $k$ whose $\bar k$-points are in $\Ga_k$-equivariant bijection with $\sA_k[2](\bar k)$; see for example \cite[Section~(3.7)(IV)]{Tat97}. 
Thus over $\mathcal{O}_K$, we deduce that
\[ \sA[2]^\et(\ov K)=\sA[2]^\et(K_{\nr})\cong \sA[2]^\et(\bar k)=\sA[2](\bar k),\]
which are all isomorphisms of (unramified) $\Ga_K$-modules.

\begin{lem} The connected-\'etale sequence
\[ 0\lra \sA[2]^\circ\lra \sA[2] \lra \sA[2]^\et \lra 0 \]
splits if and only if the exact sequence
\[ 0 \lra \sA[2]^\circ(\overline{K}) \lra A[2](\overline{K})\lra \sA[2](\bar k) \lra 0 \]
of\, $\Ga_K$-modules splits.
\end{lem}

\begin{proof}
Clearly if the connected-\'etale sequence splits, then the sequence of $\Ga_K$-module splits. Conversely, if the sequence of $\Ga_K$-modules splits, we obtain a finite flat group scheme $G\subset A[2]$ of order $|\sA[2]^\et|$ as the image of the splitting. Taking the closure of $G$ inside $\sA[2]$ then gives a finite flat group scheme $\mathcal{G}\subset \sA[2]$ of order $|\sA[2]^\et|$ mapping isomorphically onto $\sA[2]^\et$. This gives a splitting of the connected-\'etale sequence. 
\end{proof}

Let $\sA^\vee/\O_K$ be the dual abelian surface, which exists by
\cite[Theorem 8.5]{BLR}. By \cite[Corollary~1.3]{Oda69}, 
the finite flat group $k$-schemes $\sA[2]$ and $\sA^\vee[2]$ are Cartier
duals of each other; that is, $\sA[2]=\Hom_{\mathcal{O}_K}(\sA^\vee[2],\G_{m})=\Hom_{\mathcal{O}_K}(\sA^\vee[2],\mu_2)$.
If $\sA$ has {\em ordinary} reduction, that is, $\sA[2](\bar k)\simeq(\Z/2)^2$, then so does $\sA^\vee$. Thus
each of the group schemes $\sA[2]^\circ$, $\sA[2]^\et$, $\sA^\vee[2]^\circ$, $\sA^\vee[2]^\et$
has order $4$. The Cartier dual of a
commutative finite
\'etale group $\O_K$-scheme annihilated by $2$ is connected;
hence $\Hom_{\mathcal{O}_K}(\sA^\vee[2]^\et,\mu_{2})$ is naturally embedded into $\sA[2]^\circ$.
Since the orders are equal, we obtain a canonical isomorphism
\begin{equation*}
\sA[2]^\circ\cong  \Hom_{\mathcal{O}_K}(\sA^\vee[2]^\et,\mu_{2}).  
\end{equation*}
Similarly, we have $\sA^\vee[2]^\circ\cong\Hom_{\O_K}(\sA[2]^\et,\mu_{2})$, and we also have a canonical isomorphism of $\F_2[\Ga_k]$-modules
\begin{equation*}
\Hom(\sA^\vee[2](\bar k),\mu_2(\overline{K}) )\cong \sA^\vee[2](\bar k) 
\end{equation*}
because $\mu_2(\overline{K})=\F_2$ and $\wedge^2(\sA^\vee[2](\bar k))=\F_2$. Thus the connected-\'etale sequence gives rise to
an exact sequence of $\Ga_K$-modules
\begin{equation*}
0\lra \sA^\vee[2](\bar k) \lra \sA[2](\ov K)\lra \sA[2](\bar k)\lra 0.
\end{equation*}
In particular, we have the following statement.

\begin{lem} Suppose that $\sA/\mathcal{O}_K$ is a relative abelian surface with ordinary special fibre $\sA_k$. Then $\sA[2]^\circ(\overline{K})$ is unramified as a $\Ga_K$-module.
\end{lem}

\subsection{Galois action over the residue field}

If $\sA_k$ is not supersingular, then the minimal resolution $\Kum(\sA_k)\rightarrow \sA_k/G$ is a K3 surface by \cite[Theorem B]{Kat78}. There are sixteen (geometric) exceptional curves, which appear in four disjoint ${\rm D_4}$ configurations when $\sA_k$ is ordinary, and two disjoint ${\rm D}_8$ configurations when $\sA_k$ is almost ordinary. The results of Sections~\ref{sec: geom ord} and~\ref{sec: geom almost} then describe the Galois action on these curves: if $\sA_k$ is ordinary, the exceptional curves are indexed Galois-equivariantly by $\sA^\vee[2](\bar k)\times \sA[2](\bar k)$, and if $\sA_k$ is almost ordinary, they are all $k$-rational. 

\bde
Suppose that a group $\Ga$ acts on a finite set $X$.
For a prime $\ell$, we write $\Q_\ell^{X}$ for the associated $\Q_\ell$-valued permutation representation of $\Ga$.
\ede

In particular, if $X$ is the set underlying an $n$-dimensional $\F_2$-representation of $\Ga$, we have $\dim_{\Q_\ell}(\Q_\ell^{X})=2^n$. Thanks to the results of Sections~\ref{sec: geom ord} and~\ref{sec: geom almost}, we can then describe the $\Ga_k$-representation in $\H^2_\et(\Kum(\sA_k)_{\bar k},\Q_\ell)$ explicitly as follows.

\bco
Suppose that $\sA_k$ is ordinary. 
Let $V=V_1\oplus V_2$, where $V_1=\sA^\vee[2](\bar{k})$ and $V_2= \sA[2](\bar{k})$. Then for any odd prime $\ell$, there is an isomorphism
\begin{equation*}
\H^2_\et(\Kum(\sA_k)_{\bar k},\Q_\ell) \cong \H^2_\et(\sA_{\bar{k}},\Q_\ell) \oplus \Q_\ell^{V}(-1)
\end{equation*}
of\, $\Ga_k$-representations.
\eco

\begin{cor} \label{cor: kummer 2-rank 1} 
Suppose that $\sA_k$ is almost ordinary.
Then for any odd prime $\ell$, there is an isomorphism 
\[ \H^2_\et(\Kum(\sA_k)_{\bar k},\Q_\ell)\cong \H^2(\sA_{\bar k},\Q_\ell)\oplus \Q_\ell(-1)^{\oplus 16}\]
of\, $\Ga_k$-representations.
\end{cor}

\subsection{Galois action over the fraction field}

We now consider the $\Ga_K$-representation in the cohomology group $\H^2_\et(\Kum(A)_{\ov{K}},\Q_\ell)$ of the generic fibre. In this case, we have 
an isomorphism of $\Ga_K$-modules
\begin{equation*} 
\H^2_\et(\Kum(A)_{\overline{K}},\Q_\ell) \cong \H^2_\et(A_{\overline{K}},\Q_\ell) \oplus \Q_\ell^{W}(-1),
\end{equation*}
where $W=A[2](\ov{K})$ is the $\Ga_K$-representation on the $2$-torsion points of $A$. A detailed explanation of this isomorphism can be found in \cite[Lemma 4.1]{Ove21}.

We will want to compare the Galois representations in $\H^2_\et(\Kum(\sA_{\bar k}),\Q_\ell)$ and $\H^2_\et(\Kum(A)_{\ov{K}},\Q_\ell)$, at least when $\sA_k$ is non-supersingular. In the almost ordinary case, this is entirely straightforward, so we will concentrate on the case where $\sA_k$ is ordinary.

In this case, if we let $V_1=\mathcal{A}^\vee[2](\bar{k})$ and $V_2=\mathcal{A}[2](\bar{k})$, then recall from Section~\ref{sec: prelim} that we have a short exact sequence of $\F_2[\Ga_K]$-modules
\begin{equation}
0 \lra V_1  \lra W \lra V_2 \lra 0. \label{*}
\end{equation}
Let $P\subset \GL(W)$ be the parabolic subgroup
which leaves invariant the $\F_2$-subspace $V_1 \subset W$.
Let $G_W$ (respectively, $G_V$) be the image of the action of $\Ga_K$ on $W$ (respectively, on
$V=V_1\oplus V_2$). The natural projection $\pi\colon P\to \GL(V_1)\times\GL(V_2)$
induces a surjective homomorphism $G_W\to G_V$.
Since $G_V$ and $G_W$ are finite, we have $G_W\cong G_V$ if and only if $|G_W|=|G_V|$.

\bthe \label{t0}
{\samepage Suppose that $\sA_k$ is ordinary. Then  the following properties are equivalent:
\begin{enumerate}
\item \label{num: t0 1}  There is an isomorphism of\, $\Ga_K$-representations $\Q_\ell^{V}\simeq \Q_\ell^{W}$.
\item \label{num: t0 2}  The extension of\, $\F_2[\Ga_K]$-modules $(\ref{*})$ is split.
\item \label{num: t0 3} $|G_W|=|G_V|$.
\end{enumerate}
\noindent If these properties hold, then the $\Ga_K$-representation in $W=A[2](\ov{K})$ is unramified; that is, it factors through $\Ga_k$.
}
\ethe

\begin{proof} The implications~\eqref{num: t0 1} $\Rightarrow$~\eqref{num: t0 3} and~\eqref{num: t0 2} $\Rightarrow$~\eqref{num: t0 1} are clear.
Either~\eqref{num: t0 1} or~\eqref{num: t0 2} implies that the representation of $\Ga_K$ in $W$ is unramified because
the representation of $\Ga_K$ in $V$ is unramified. It remains to show that~\eqref{num: t0 3} implies~\eqref{num: t0 2}.

Choosing a section $V_2\rightarrow W$ of the projection $W\to V_2$
induces a section $\sigma$ of the
projection $\pi\colon P\to \GL(V_1)\times\GL(V_2)$. Thus we have a split exact sequence
of groups
\begin{equation}
0\lra R_u(P)\lra P\lra  \GL(V_1)\times\GL(V_2)\lra 0,\label{ext}
\end{equation}
where $R_u(P)$ is the unipotent radical of $P$. 
We need to show that if
$G_W$ is isomorphic to $G_V=\pi(G_W)$, then $G_W$ is conjugate, via an element of $R_u(P)$, to 
$\sigma(G_V)$. 

The pullback of \eqref{ext} with respect to the inclusion $G_V\subset  \GL(V_1)\times\GL(V_2)$ is a split exact sequence
\[0\lra R_u(P)\lra R_u(P)\rtimes G_V\lra G_V\lra 0,\]
where $G_V\subset \GL(V_1)\times\GL(V_2)$ acts on $R_u(P)$
by conjugation.
The subgroups $G_W$ and $\si(G_V)$ are complements to $R_u(P)$ in $R_u(P)\rtimes G_V$.
It is enough to show that $\H^1(G_V,R_u(P))=0$ because 
then all the complements to $R_u(P)$ in $R_u(P)\rtimes G_V$ are conjugate; see, \textit{e.g.},~\cite[Proposition 9.21]{Rot09}. 

We note that the action of $(A,B)\in\GL(V_1)\times\GL(V_2)$ on $\left(\begin{smallmatrix}
I & X \\ 0 & I
\end{smallmatrix}\right)\in R_u(P)$ by conjugation
via $\si((A,B))$ sends $X$ to $AXB^{-1}$, so $R_u(P)$ is isomorphic as a $\GL(V_1)\times\GL(V_2)$-representation to $V_1\otimes_{\F_2} V_2^*\cong V_1\otimes_{\F_2} V_2$.
It remains to prove that for any subgroup $G\subset \GL(V_1)\times\GL(V_2)$, we have 
$\H^1(G,V_1\otimes_{\F_2} V_2)=0$. Since this group is annihilated by $2$, by the standard restriction-corestriction
argument, it is enough to prove that $\H^1(H,V_1\otimes_{\F_2} V_2)=0$, 
where $H$ is a $2$-Sylow subgroup of $G$.
Then $H$ is contained in a subgroup of $S_3\times S_3$ isomorphic to the product of 
$\Z/2\subset S_3$ and $\Z/2\subset S_3$.
The $\F_2$-vector space $V_1$ has a basis whose elements are permuted by $\Z/2$, and similarly for $V_2$.
This gives a basis of $V_1\otimes_{\F_2} V_2$ whose elements are permuted by
$\Z/2\times \Z/2$. It follows that $\H^1(H,V_1\otimes_{\F_2} V_2)=0$ 
for any subgroup $H\subset\Z/2\times \Z/2$.
\end{proof}

\section{Explicit smooth models} \label{sec: explicit}

Let $\sA/\O_K$ be an abelian scheme of relative dimension $2$ whose reduction $\sA_k$ is non-supersingular. Let $A:=\sA_K$ be the generic fibre of $\sA$. In this section, we show how to construct an explicit smooth model for $\Kum(A)$ under suitable assumptions on the $2$-torsion $A[2](\overline{K})$ as a $\Ga_K$-module. In particular, this will result in a proof of Theorem~\ref{theo: main 1}.

\subsection{Blowups and specialisation}

Let $\mathcal{X}\rightarrow \Spec(\O_K)$ be a flat morphism of finite type, of relative dimension $2$, and with normal, integral fibres. It is not true in general that blowing up a closed subscheme of $\mathcal{X}$ commutes with base change to $k$, even if the centre is flat over $\O_K$. However, we do have the following result. 

\bpr \label{lem: rdp blowups}
Suppose that $\mathcal{Z}\subset \mathcal{X}$ is an $\O_K$-section such that both fibres of $\mathcal{X}$ have an isolated rational double point at $\mathcal{Z}$. Then $\Bl_{\mathcal{Z}}(\mathcal{X})_k\cong \Bl_{\mathcal{Z}_k}(\mathcal{X}_k)$.
\epr

\begin{proof}
Let $\mathcal{I}\subset \O_\mathcal{X}$ be the ideal of $\mathcal{Z}$ and $\mathcal{I}_k\subset \mathcal{O}_{\mathcal{X}_k}$ the ideal of $\mathcal{Z}_k$ in $\mathcal{X}_k$. Thus $\mathcal{I}_k$ is just the image of $\mathcal{I}$ inside $\O_{\mathcal{X}_k}$, and similarly each power $\mathcal{I}_k^n$ is the image of the corresponding power $\mathcal{I}^n$ inside $\O_{\X_k}$. Then
\[ \Bl_{\mathcal{Z}}(\mathcal{X})_k = \Proj\left( \bigoplus_{n\geq0} \mathcal{I}^n \otimes_{\O_K} k \right), \] 
whereas
\[ \Bl_{\mathcal{Z}_k}(\mathcal{X}_k) = \Proj\left( \bigoplus_{n\geq0} \mathcal{I}_k^n  \right),  \]
so we need to show that the natural map
$\mathcal{I}^n \otimes_{\O_K} k \longrightarrow \mathcal{I}_k^n$
is an isomorphism for all $n$. Since this map is clearly surjective, we see that the natural map
\[ \mathcal{O}_{\mathcal{X}}/\mathcal{I}^n \otimes_{\O_K} k \to \O_{\X_k}/\mathcal{I}^n_k \] 
is an isomorphism. Hence showing that $\mathcal{I}^n \otimes_{\O_K} k \rightarrow \mathcal{I}_k^n$ is an isomorphism is equivalent to showing that the exact sequence
\[ 0 \lra \mathcal{I}^n\lra \O_\mathcal{X} \lra \mathcal{O}_{\mathcal{X}}/\mathcal{I}^n \lra 0 \]
remains exact upon applying $-\otimes_{\O_K}k$, which will certainly follow if we can show that $\mathcal{O}_{\mathcal{X}}/\mathcal{I}^n$ is flat over $\O_K$. 

To see this, we first observe that each $\mathcal{O}_{\mathcal{X}}/\mathcal{I}^n$ is finite over $\O_K$. Indeed, $\mathcal{O}_{\X}/\mathcal{I}$ is finite over $\O_K$, and $\mathcal{I}/\mathcal{I}^2$ is a finitely generated $\mathcal{O}_\X/\mathcal{I}$-module. Hence the surjective map $(\mathcal{I}/\mathcal{I}^2)^{\otimes n}\to \mathcal{I}^n/\mathcal{I}^{n+1}$ shows that each $\mathcal{I}^n/\mathcal{I}^{n+1}$ is a finitely generated $\O_K$-module, and thus so is each $\mathcal{O}_{\mathcal{X}}/\mathcal{I}^n$ by induction on $n$.

By the structure theorem for modules over a principal ideal domain, we see that each $\mathcal{O}_{\mathcal{X}}/\mathcal{I}^n$ is isomorphic to an $\O_K$-module of the form
\[ \O_K^{\oplus r}\oplus \O_K/\mathfrak{m}_K^{n_1} \oplus \dots \oplus \O_K/\mathfrak{m}_K^{n_s}. \]
Thus flatness over $\O_K$ is equivalent to having $s=0$, which in turn is equivalent to the equality of dimensions
\begin{equation} \label{eqn: dims} \dim_K \mathcal{O}_{\mathcal{X}_K}/\mathcal{I}_K^n=\dim_k \mathcal{O}_{\mathcal{X}_k}/\mathcal{I}_k^n. \end{equation}
Since $\mathcal{I}_K\subset \mathcal{O}_{\mathcal{X}_K}$ and $\mathcal{I}_k\subset \mathcal{O}_{\mathcal{X}_k}$ are both maximal ideals defining isolated rational double points on normal surfaces, it follows from Lemma~\ref{lem: rdp dim} below that both sides of \eqref{eqn: dims} are equal to $n^2$.
\end{proof}

\begin{lem} \label{lem: rdp dim} Let $(A,\mathfrak{m}_A)$ be a complete Noetherian local ring, normal of Krull dimension $2$, with a rational double point at $\mathfrak{m}_A$. Then $\dim_{A/\mathfrak{m}_A} (\mathfrak{m}_A^n/\mathfrak{m}_A^{n+1})= 2n+1$. 
\end{lem}

\begin{proof}
We may assume that the residue field $A/\mathfrak{m}_A$ is algebraically closed. In this case, the computation is done in Theorem 4 and Corollary 6 of \cite{Art66}.
\end{proof}

This result has the following important consequence. 

\begin{lem} \label{lemma: key!} Let $x\in \mathcal{X}(k)$, and assume that $\mathcal{X}_k$ has a rational double point of type ${\rm A}_1$ at $x$. Then there exists at most one $($geometric$)$ rational double point of $\mathcal{X}_K$ specialising to $x$. 
\end{lem}

\begin{rem} 
If $P\colon \Spec(\overline{K})\to \X_K$ is a point, then there exist a finite extension $L/K$ and a finitely generated $\O_K$-algebra $R\subset L$ such that $R[1/2]=L$ and the scheme-theoretic image of $P$ is given by a closed immersion $\Spec(R) \to \X$.
Since the normalisation of $R$ is equal to $\O_L$, the special fibre $\Spec(R\otimes_{\O_K} k)$ has a single point. We say that $P$ specialises to $x$ if the closed immersion $\Spec(R)\to \X$ sends the unique point of $\Spec(R\otimes_{\O_K} k)$ to $x$.
\end{rem}

\begin{proof}
Suppose towards a contradiction that there are two distinct such points $P_1$, $P_2$. Let $L/K$ be a finite extension over which both are defined, with ring of integers $\O_L$ and residue field $k_L$. Let $\mathcal{X}_{\O_L}$ be the base change of $\mathcal{X}$ to $\mathcal{O}_L$. Thus $\mathcal{X}_{\O_L}$ is flat over $\mathcal{O}_L$ with normal fibres, its special fibre $\mathcal{X}_{k_L}$ has a rational double point of type ${\rm A}_1$ at $x$, and both $P_1$ and $P_2$ specialise to this point. 

Let $\mathcal{P}_1$ denote the scheme-theoretic closure of $P_1$ inside $\mathcal{X}_{\O_L}$. This is therefore an $\mathcal{O}_L$-section of $\mathcal{X}_{\O_L}$, and we let $\mathcal{X}_{\O_L}'$ denote the blowup of $\mathcal{X}'$ along this section. Thus the special fibre $\mathcal{X}_{k_L}'$ is smooth over $k_L$ at every point in the fibre over $x$. Now, we let $\mathcal{P}_2$ denote the closure of $P_2$ inside $\mathcal{X}_{\O_L}$ and $\mathcal{P}_2'\subset \mathcal{X}_{\O_L}'$ its strict transform. Then $\mathcal{P}_2'$ is an $\mathcal{O}_L$-section of $\mathcal{X}_{\O_L}'$ such that the special fibre $\mathcal{X}_{k_L}'$ is smooth at $\mathcal{P}'_{2,k_L}$, but the generic fibre $\mathcal{X}'_{L}$ is singular at $\mathcal{P}'_{2,L}$. This is the contradiction we seek.
\end{proof}

\begin{rem} More generally, if we have a collection of rational double points, we can talk about their `total degree' as being the sum of the lower indices in the ADE classification. We can then use Lemma~\ref{lemma: key!} to prove by induction that the total degree of a collection of  rational double points cannot increase under specialisation. We will not need this more general result.
\end{rem}

\subsection{Blowing up $\bm{\sA/G}$ along the \'etale part of $\bm{2}$-torsion}

We now let $\sY:=\sA/G$ be the quotient scheme by the involution $\iota$ and let
$q\colon\sA\rightarrow \sY$ be the quotient morphism. Note that $\mathcal{Y}$ is flat over $\O_K$.
Our aim is to show that under appropriate conditions one can explicitly construct 
relative (smooth) Kummer surfaces by resolving the relative quotient surface $\mathcal{Y}$.

We begin by describing the fibres of $\sY\to\Spec(\O_K)$.
Since quotients of quasi-projective schemes by finite group actions commute with flat base change,\footnote{This is essentially because the ring of functions on the quotient can be expressed as a kernel; see for example the argument in Section~(4.24) on p.~59 of \cite{EvdGM}.}
we get a natural  identification $A/G\, \tilde\to\, \mathcal{Y}_K $. In fact, the same is true on the special fibre. 

\begin{prop} \label{spec quot}
The natural map $\sA_k\rightarrow \mathcal{Y}_k$ induces an isomorphism $
\sA_k/G\,\tilde\to\,\mathcal{Y}_k$.
\end{prop}

\begin{proof}
In view of compatibility of quotients by $G$ with flat base change, we may assume that $k$ is algebraically closed. Since the $G$-action is free on $\sA\setminus \sA[2]$, we obtain an isomorphism
\[ (\sA_k\setminus \sA_k[2])/G \,\tilde\lra\, (\sA/G)_k \setminus q(\sA_k[2]).  \]
Once more appealing to compatibility with flat base change, it therefore suffices to prove the analogous statement after base changing to the formal completion $(\sA/G)_{/Q}$ at any point $Q$ of $q(\sA_k[2])$. Up to translation, and possibly enlarging $K$ if necessary to ensure that the map $\sA[2](K)\rightarrow \sA[2](k)$ is surjective, we may therefore replace $\sA$ with its completion at the zero section $O_k\in \sA_k$ on the special fibre.

Thus we have an action of $G$ on $\widehat{\mathcal{O}}_{\sA,O_k}\cong \mathcal{O}_K\pow{u,v}$, and we may choose the local parameters $u$ and $v$ in such a way that the description of the action and its quotient given in \cite[Proposition 1.1]{Sch09} and \cite[Theorem, p.~60]{Art75} holds modulo $\mathfrak{m}_K$. The quotient $\mathcal{O}_K\pow{u,v}/ \O_K\pow{u,v}^G$ is $\mathfrak{m}_K$-torsion free, which shows that the map $\O_K\pow{u,v}^G\to \O_K\pow{u,v}$ remains injective upon applying $-\otimes_{\O_K}k$. This implies that the natural map $\mathcal{O}_K\pow{u,v}^G \otimes_{\mathcal{O}_K} k\rightarrow k\pow{u,v}^G$ is injective; we need to show that it is in fact an isomorphism. In other words, we need to show that
\[ \mathcal{O}_K\pow{u,v}^G\lra k\pow{u,v}^G\]
is surjective. But this is straightforward since every element of $k\pow{u,v}^G$ can be written as a series in $x=u\cdot \iota(u)$, $y=v\cdot \iota(v)$, and $z=u\cdot \iota(v)+v\cdot \iota(u)$. Any such series clearly lifts to an element of $\mathcal{O}_K\pow{u,v}^G$.
\end{proof}

Let us now assume that $\sA_k$ is non-supersingular and that the exact sequence of $\Ga_K$-modules
\begin{equation*}
 0\lra \sA[2]^\circ(\ov K) \lra \sA[2](\overline{K}) \lra \sA[2](\bar k) \lra 0
\end{equation*}
is split. Then the surjective morphism of finite flat group $\O_K$-schemes
$\sA[2]\to\sA[2]^\et$ has a section $\si$, so that $\sA[2]$ is isomorphic to 
$\sA[2]^\circ\times \si(\sA[2]^\et)$. The degree of $\si(\sA[2]^\et)\to\Spec(\O_K)$ is $4$ in the ordinary case and $2$ in the almost ordinary case. Let $\mathcal{Z}=q(\si(\sA[2]^\et))\subset \mathcal{Y}$ denote the scheme-theoretic image of $\si(\sA[2]^\et)$.

\begin{lem} 
The natural map $\si(\sA[2]^\et)\rightarrow \mathcal{Z}$ is an isomorphism. 
\end{lem}

\begin{proof}
The claim can be checked after making a finite extension of $K$, so we may assume that $\si(\sA[2]^\et) \subset \sA$ consists of either four (in the ordinary case) or two (in the almost ordinary case) disjoint $\mathcal{O}_K$-sections. Their images remain disjoint in $\mathcal{Y}$, and so we may reduce to considering the scheme-theoretic image in $\mathcal{Y}$ of a single $\mathcal{O}_K$-section in $\sA$, where the claim is clear.
\end{proof} 

Let $\mathcal{Y}^{(1)}$ be the blowup of $\mathcal{Y}$ in $\mathcal{Z}$. It then follows from 
Proposition~\ref{lem: rdp blowups} that the special fibre $\mathcal{Y}^{(1)}_k$ is the blowup of $\mathcal{Y}_k$ in its reduced singular locus.

\subsection{The ordinary case}

We now assume further that $\sA_k$ is ordinary. Thus both the generic fibre $\mathcal{Y}^{(1)}_K$ and the special fibre $\mathcal{Y}^{(1)}_k$ contain twelve (geometric) singular points, all of which are rational double points of type ${\rm A}_1$. Each of the twelve singular points on $\mathcal{Y}^{(1)}_K$ has to specialise to a singular point on $\mathcal{Y}^{(1)}_k$, and by Lemma~\ref{lemma: key!} each singular point on $\mathcal{Y}^{(1)}_k$ is specialised to by at most one singular point on $\mathcal{Y}^{(1)}_K$. 

In other words, if we let $\mathcal{Z}^{(1)} \subset \mathcal{Y}^{(1)}$ denote the scheme-theoretic closure of the twelve (geometric) singular points on $\mathcal{Y}^{(1)}_K$, then $\mathcal{Z}^{(1)}$ is finite \'etale over $\O_K$ of degree $12$ and intersects both the generic and special fibres in their reduced singular loci. We now let $\X$ be the blowup of $\mathcal{Y}^{(1)}$ in $\sZ^{(1)}$. 

\bthe \label{ord}
Suppose that $\sA_k$ is ordinary and the exact sequence
\[0\lra \sA[2]^\vee(\bar k) \lra \sA[2](\overline{K}) \lra \sA[2](\bar k) \lra 0\]
splits. Then any section $\sigma$ of $\sA[2]\to\sA[2]^\et$ gives rise to
a smooth and projective scheme $\X/\O_K$, equipped with an action of $\sA[2]^\et$, whose fibres are the minimal desingularisations of the fibres of $\sA/G$. In particular, 
$\X_K\cong\Kum(\sA_K)$ and $\X_k\cong\Kum(\sA_k)$.
\ethe
\begin{proof} 
The statement about the special fibre follows from the results of Section~\ref{sec: geom ord}; in particular, $\X_k$ is smooth. The generic fibre $\X_K$ is clearly the Kummer surface attached to $\sA_K$, so it is smooth. We can see that $\X/\O_K$ is flat because $\X$ is integral, $\O_K$ is a DVR, and $\X\rightarrow \Spec(\O_K)$ is surjective; see \cite[Proposition~III.9.7]{Har77}. It is projective because it is an iterated blowup of the projective $\O_K$-scheme $\mathcal{Y}$.
Finally, $\X/\O_K$ is smooth since it is flat with smooth fibres; see
\cite[\href{https://stacks.math.columbia.edu/tag/01V8}{Lemma 01V8}]{stacks-project}. 
\end{proof}

\subsection{The almost ordinary case}

The construction of an explicit smooth resolution of $\mathcal{Y}=\sA/G$ is slightly more involved in the almost ordinary case. Here, we assume that $\sA[2](\overline{K})=\sA[2](K)$, and we show that we can perform the `na\"ive' resolution described in Section~\ref{sec: naive} on the relative surface $\mathcal{Y}$. In fact, we will replace $\mathcal{A}$ with the formal completion $\mathcal{S}= \Spec(\widehat{\mathcal{O}}_{\sA,O_k})$ at the zero section of the special fibre. We let $\T=\mathcal{S}/G$ and write $q\colon \mathcal{S}\rightarrow \T$ for the quotient map. We will describe a resolution of $\mathcal{S}$ via an explicit sequence of blowups. To obtain a resolution of $\mathcal{Y}$, we simply perform the same sequence of blowups and then translate the whole procedure by the point of $\sA[2](K)$ given by the image of the non-identity point of $\sA[2](k)$ under the splitting $\sigma$ of $\sA[2]\to\sA[2]^\et$. 

First of all, note that any point in the kernel of the reduction map $\sA(K)\rightarrow \sA(k)$ can naturally be thought of as a $K$-point of $\mathcal{S}$. Moreover, any such point extends uniquely to a section $\Spec(\O_K)\rightarrow \mathcal{S}$. In particular, all eight points of $\sA[2]^\circ(K)$ give rise to sections of $\mathcal{S}$, and the images of these sections under $q$ intersect the generic fibre $\T_K$ precisely at its eight singular points. The fibre product $\mathcal{T}^{(1)}:=\mathcal{Y}^{(1)}\times_\mathcal{Y} \T$ is then the blowup of $\T$ along $q(O)$.

We therefore know that the (reduced) singular locus of $\T^{(1)}_K$ consists of seven $K$-rational points of type ${\rm A}_1$, and the (reduced) singular locus of $\T^{(1)}_k$ consists of two $k$-rational points, one of type ${\rm A}_1$ and one of type ${\rm D}_6^1$. By Lemma~\ref{lemma: key!}, there exists at least one singular point on $\T^{(1)}_K$ specialising to the ${\rm D}_6^1$-singularity $Q_2^{(1)}$ (in the notation of Section~\ref{sec: geom almost}). Let $\mathcal{Q}_2^{(1)}\subset \T^{(1)}$ denote the closure of this point and $\mathcal{T}^{(2)}$ the blowup of $\T^{(1)}$ along $\mathcal{Q}_2^{(1)}$. 

Now the (reduced) singular locus of $\T^{(2)}_K$ consists of six $K$-rational points of type ${\rm A}_1$, and the (reduced) singular locus of $\T^{(2)}_k$ consists of three $k$-rational points, two of type ${\rm A}_1$ and one of type ${\rm D}_4^0$. By Lemma~\ref{lemma: key!}, there exists at least one singular point on $\T^{(2)}_K$ specialising to the ${\rm D}_4^0$-singularity $Q_3^{(2)}$. We let $\mathcal{Q}_3^{(2)}\subset \T^{(2)}$ denote the closure of this point and $\mathcal{T}^{(3)}$ the blowup of $\T^{(2)}$ along $\mathcal{Q}_3^{(2)}$. 

We therefore see that both the general and special fibres of $\T^{(3)}$ have precisely five singular points, all of type ${\rm A}_1$, and the singular subscheme of $\T^{(3)}$ is finite \'etale of degree $5$ over $\O_K$.
Now, as in the ordinary case, we may blow up the singular subscheme of $\T^{(3)}$ to obtain a scheme which is formally smooth over $\O_K$. We therefore obtain the following analogue of Theorem~\ref{ord}.

\bthe \label{alm}
Suppose that $\sA_k$ is almost ordinary and that $\sA[2](\overline{K})$ is trivial as a $\Ga_K$-module. Then there exists a smooth and projective scheme $\X/\O_K$, equipped with an action of $\sA[2]^\et$, whose fibres are the minimal desingularisations of the fibres of $\sA/G$.  In particular, 
$\X_K\cong\Kum(\sA_K)$ and $\X_k\cong\Kum(\sA_k)$.
\ethe

\brem It is not completely transparent exactly where the hypothesis that $\sA[2](\overline{K})=\sA[2](K)$ is used. However, a more careful examination of the proof shows that:
\begin{enumerate}
\item the two sections $\mathcal{Q}_2^{(1)}$ and $\mathcal{Q}_3^{(2)}$ that we blew up above,
\item the chosen splitting of $\sA[2]\to \sA[2]^\et$,
\item the canonical multiplicative subgroup $\mu_2\cong \Hom_{\O_K}(\sA^\vee[2]^\et,\mu_2) \subset \sA[2]$
\end{enumerate}
all combine to give a set of $K$-rational generators of $\sA[2](\overline{K})$. The choice involved in producing $\mathcal{X}$ is essentially a choice of $\F_2$-basis of $\sA[2](\overline{K})$ containing the non-trivial element of $\Hom_{\O_K}(\sA^\vee[2]^\et,\mu_2)$. 
\erem

\section{Good reduction criterion}
 \label{sec: grc}
 
In this section, we complete the proofs of Theorems~\ref{theo: main 1},~\ref{theo: main 2}, and~\ref{theo: main 3}. Since we will want to use the main result of \cite{CLL}, we first need to explain the `canonical reduction' of a K3 surface which plays a key role there. 

\subsection{Canonical reduction}

Suppose that $X/K$ is a K3 surface which attains good reduction after a finite and unramified extension of~$K$. In general, this does not imply that $X$ has good reduction over $K$ itself; see for example \cite[Section~7]{LM18}. However, it is still possible to produce a K3 surface $X_0/k$ as the `reduction' of $X$, in a way that is unique up to $k$-isomorphism. 

The key result that allows us to do this is \cite[Proposition 4.7(2)]{LM18}, which says that if $\mathcal{X}_1$ and $\mathcal{X}_2$ are smooth models of our K3 surface $X$ over $\O_K$, that is,
smooth and proper algebraic spaces over $\O_K$ with generic fibres identified with $X$, then $\mathcal{X}_1$ and $\mathcal{X}_2$ are connected by a sequence of flopping contractions and their inverses. In particular, this implies that the canonical rational map $\mathcal{X}_1\dashrightarrow \mathcal{X}_2$ given by the identity on generic fibres is an isomorphism away from a finite collection of curves on the special fibres of $\mathcal{X}_1$ and $\mathcal{X}_2$. We may therefore restrict this map to obtain a birational map $\mathcal{X}_{1,k}\dashrightarrow \mathcal{X}_{2,k}$ between these special fibres. As remarked in the introduction, these special fibres must be K3 surfaces (since they have trivial canonical bundle and coherent Euler characteristic $2$), and hence this restricted birational map is in fact an isomorphism.

This has the following important consequence. Suppose that $L/K$ is a finite and unramified Galois extension, with induced residue field extension $k_L/k$. If $\mathcal{X}/\O_L$ is a smooth model for $X_L$, then we may consider the natural $\Gal(L/K)$-action on $X_L$ as a rational action on $\mathcal{X}$. By the above discussion, this is defined away from a finite collection of curves on $\mathcal{X}_{k_L}$. We can therefore restrict this rational action to the special fibre $\mathcal{X}_{k_L}$, and again, as in the above discussion, this restricted rational action is regular. Hence we may form the quotient $X_0:=\mathcal{X}_{k_L}/\Gal(L/K)$, which is a K3 surface over $k$ by the theory of Galois descent. Moreover, if we had any two such smooth models $\mathcal{X}_1$ and $\mathcal{X}_2$ over $\O_L$, then the identity map between the generic fibres induces a $\Gal(L/K)$-equivariant isomorphism between their special fibres. Thus the K3 surface $X_0$, up to isomorphism, does not depend on the choice of model $\mathcal{X}/\O_L$, or indeed on the choice of $L$.

\begin{defi} The K3 surface $X_0$ over $k$ is called the {\em canonical reduction} of $X$.
\end{defi}

In our case the canonical reduction appears in the following way.

\begin{lem}\label{lem: can red alt}  
Let $\mathcal{Y}/\O_K$ be a flat, projective scheme with normal,
$2$-dimensional fibres, such that the minimal resolution $X$ of\,
$\mathcal{Y}_K$ is a K3 surface.  Assume that there exist a finite and
unramified Galois extension $L/K$ and a proper birational morphism
$\X\to\mathcal{Y}_{O_L}$ such that $\X$ is smooth over $\O_L$ with
generic fibre isomorphic to $X_L$.  Then the minimal resolution of\,
$\mathcal{Y}_k$ is the canonical reduction of\, $X$.
\end{lem}

\begin{proof}
The Galois group $\Gal(L/K)$ acts naturally on $X_L$ and $\mathcal{Y}_L$, and the morphism
$X_L\to \mathcal{Y}_L$ is $\Gal(L/K)$-equivariant. The action of $\Gal(L/K)$ on $X_L$ extends uniquely to a rational action of $\Gal(L/K)$ on $\X$. The morphism $\X\to\mathcal{Y}_{O_L}$ is thus $\Gal(L/K)$-equivariant. 
As we have seen, the rational action of $\Gal(L/K)$ on $\X$ restricts to a regular action on $\X_{k_L}$. We therefore obtain a $\Gal(k_L/k)$-equivariant birational morphism $\mathcal{X}_{k_L}\rightarrow  \mathcal{Y}_{k_L}$. 
Since $\mathcal{X}_{k_L}$ is a K3 surface, this morphism is a
minimal resolution of $\mathcal{Y}_{k_L}$.
Hence $\mathcal{X}_{k_L}/\Gal(L/K)$, which by definition is the canonical reduction of $X$, is the minimal resolution of $\mathcal{Y}_k$.
\end{proof}

\subsection{The proofs}

We can now prove our main theorems.

\bthe \label{t1}
Let $A$ be an abelian surface over $K$ with good, non-supersingular reduction.
Then we have the following statements:
\begin{enumerate}
\item  \label{num: t1 1} The Kummer surface $\Kum(A)$ attains good reduction
with a scheme model after a finite field extension $L/K$ of degree $[L:K]\leq |\GL_4(\F_2)|$.
\item
\label{num: t1 2}
If the $\Ga_K$-module $A[2](\ov{K})$ is unramified, then the extension $L/K$ 
in \eqref{num: t1 1} can be chosen to be unramified.
\item \label{num: t1 3}
  \begin{enumerate}[label={\rm(\alph*)}, ref={\rm 3\alph*}]
\item\label{num: t1 3 a}  If $A$ has ordinary reduction, then $\Kum(A)$ has good reduction over $K$ if and only if 
the exact sequence
\[ 0 \lra \sA^\vee[2](\bar k) \lra A[2](\ov{K})\lra \sA[2](\bar k) \lra 0 \]
of\,  $\Ga_K$-modules splits.
If this condition holds, then the $\Ga_K$-module $A[2](\ov{K})$ is unramified and $\Kum(A)$ has good reduction over $K$ with a scheme model.
\item\label{num: t1 3 b} If $A$ has almost ordinary reduction, then $\Kum(A)$ has good reduction over $K$ if and only if the $\Ga_K$-module $A[2](\ov{K})$ is trivial.
If this condition holds, then $\Kum(A)$ has good reduction over $K$ with a scheme model.
  \end{enumerate}
\end{enumerate}
\ethe

\begin{proof} Parts~\eqref{num: t1 1} and~\eqref{num: t1 2}, as well as `if' statements of part~\eqref{num: t1 3} immediately follow from Theorems~\ref{ord} and~\ref{alm}. 
  In particular, we see that the K3 surface $\Kum(A)$ satisfies Hypothesis~($\star$) in the terminology of \cite{LM18} and \cite{CLL}.

For the `only if' direction of part~\eqref{num: t1 3}, we first note that the hypothesis that $\Kum(A)$ has good reduction implies that $A[2](\overline{K})$ is an unramified $\Ga_K$-module. We then consider the singular model $\sA/\iota$ for $\Kum(A)$, which by the results of Section~\ref{sec: explicit} has a simultaneous resolution after an unramified extension of $K$. Hence Lemma~\ref{lem: can red alt} shows that $\Kum(\sA_k)$ is the canonical reduction of $\Kum(A)$.

It therefore follows from \cite[Theorem 1.6]{CLL} that good reduction of $\Kum(A)$ over $K$ is equivalent to the existence of an isomorphism of $\Ga_K$-representations
\[\H^2_\et(\Kum(A)_{\ov K},\Q_\ell)\simeq \H^2_\et(\Kum(\sA_k)_{\bar k},\Q_\ell)\]
for any odd prime $\ell$.  Indeed, given this isomorphism, we have an induced isomorphism on respective semisimplifications and therefore an isomorphism
\[ \frac{\H^2_\et(\Kum(A)_{\ov K},\Q_\ell)}{\H^2_\et(A_{\ov K},\Q_\ell)} \simeq  \frac{\H^2_\et(\Kum(\sA_k)_{\bar k},\Q_\ell)}{\H^2_\et(\sA_{\bar k},\Q_\ell)} \]
since both sides are semisimple and $\H^2_\et(A_{\ov K},\Q_\ell)\cong \H^2_\et(A_{\bar k},\Q_\ell)$. Hence we can apply Theorem~\ref{t0} in the case of ordinary reduction or Corollary~\ref{cor: kummer 2-rank 1} in the case of almost ordinary reduction.
\end{proof}

\begin{ex}
 \label{exa: nonsplit} 
The condition appearing in Theorem~\ref{t1}\eqref{num: t1 3 a} is not automatic, even if we assume that $A[2](\overline{K})$ is unramified. To see this, we let $\mathcal{E}/\Z_2$ be the elliptic curve defined by
\[\mathcal{E}: y^2+xy=x^3-20x-5. \]If we reduce modulo $2$, we get the curve
\[ \mathcal{E}_{\F_2}:y^2+xy=x^3+1, \]
which is smooth over $\F_2$, and hence $\mathcal{E}$ is indeed an elliptic curve over $\Z_2$. Note that $\mathcal{E}_{\F_2}$ has the $2$-torsion point $(x,y)=(0,1)$ and is therefore ordinary. If we make the change of co-ordinates $y=\frac{1}{2} (y'-x)$, then $E:=\mathcal{E}_{\Q_2}$ can be defined by the equation
\[ y'^2=4x^3+x^2-80x-20, \]
and the right-hand side factors as
\[4x^3+x^2-80x-20=(4x+1)(x^2-20). \]
Since $20$ is not a square in $\Q_2$, it follows that $E(\Q_2)$ has precisely one point of exact order $2$, namely $(x,y)=\left(-\frac{1}{4},\frac{1}{8}\right)$. However, the full $2$-torsion of $E$ is defined over the unramified extension $\Q_2(\sqrt{5})$ of $\Q_2$. The exact sequence
\[ \xymatrix{ 0 \ar[r] & \mathcal{E}[2]^\circ(\overline{\Q}_2) \ar@{=}[d] \ar[r] &  E[2](\overline{\Q}_2) \ar[r] & \mathcal{E}[2](\overline{\F}_2) \ar[r] & 0 \\ &  \mathcal{E}[2](\overline{\F}_2) } \] therefore consists of unramified $\F_2[\Ga_K]$-modules, since all points of $E[2]$ are defined over $\Q_2(\sqrt{5})$, and is non-split since $\mathcal{E}[2](\overline{\F}_2)$ is a trivial $\Ga_k$-module but $E[2](\overline{\Q}_2)$ is not. In particular, if we take $\mathcal{A}=\mathcal{E}\times_{\Z_2}\mathcal{E}$, then $A:=\mathcal{A}_{\Q_2}$ is an abelian surface with good ordinary reduction, and the exact sequence
\[0 \lra \mathcal{A}^\vee[2](\overline{\F}_2) \lra \mathcal{A}[2](\overline{\Q}_2) \lra \mathcal{A}[2](\overline{\F}_2) \lra 0 \]
consists of unramified $\F_2[G_K]$-modules and is non-split. Thus $\Kum(A)$ has good reduction over $\Q_2(\sqrt{5})$ but not over $\Q_2$.
\end{ex}

We leave it to the reader to find an example of an abelian surface $A/K$ with good, almost ordinary reduction, such that $A[2](\overline{K})$ is unramified but non-trivial as a $\Ga_K$-module.

\begin{rem} \label{height}
In the case where $\Kum(A)$ has good reduction, we can ask what the {\em height} is of the reduced K3 surface over $k$, that is, the height of its formal Brauer group over $k$ as defined in \cite{AM77}. This is the same as the height of the Kummer surface $\Kum(\sA_k)$ and can be detected in the slopes (roughly speaking, the valuations of the Frobenius eigenvalues) of the geometric crystalline cohomology groups of $\Kum(\sA_k)$. After replacing $k$ with $\bar{k}$, and using crystalline cohomology relative to $W:=W(\bar k)$, we have an isomorphism of $F$-isocrystals
\[ \H^2_{\cris}(\Kum(\sA_{\bar k})/W)_{\Q} \cong \H^2_{\cris}(\sA_{\bar k}/W)_{\Q} \oplus W(-1)_{\Q}^{\oplus 16}.   \]
This allows us to simply read off the slopes  of $\H^2_{\cris}(\Kum(\sA_{\bar k})/W)_{\Q}$ from those of
\[ \H^2_{\cris}(\sA_{\bar k}/W)_{\Q}= \bigwedge^2 \H^1_{\cris}(\sA_{\bar k}/W)_{\Q}.\]
Indeed, if $\sA_k$ has $2$-rank $r\in \{1,2\}$, then $\H^1_{\cris}(\sA_{\bar k}/W)_{\Q}$ has slopes $0$, $\frac{1}{2}$, $1$ with multiplicities $r$, $4-2r$, $r$, respectively. If we therefore set $h:=3-r$, then $\H^2_{\cris}(\sA_{\bar k}/W)_{\Q}$ has slopes $1-\frac{1}{h}$, $1$, $1+\frac{1}{h}$ with multiplicities $h$, $6-2h$, $h$, respectively. Since $W(-1)_{\Q}$ has slope $1$, we deduce from \cite[Section~II.7.2]{Ill79} that $\Kum(\sA_k)$ has height~$1$ when $\sA_k$ is ordinary, and height $2$ when $\sA_k$ is almost ordinary. 
\end{rem}

\section{Twisted Kummer surfaces}

We can also play the same game with Kummer surfaces obtained via twisting. Indeed, let $A/K$ be an abelian surface, and let $Z$ be a $K$-torsor for $A[2]$. 
We write $[Z]$ for the class of $Z$ in $\H^1(K,A[2])$.
The quotient $A_Z:= (A\times_K Z)/A[2]$, where $A[2]$ acts diagonally on the product, is a $K$-torsor for $A$. There is an action of $G$ on $A_Z$, coming from the $G$-action on $A$, and we can form the quotient $A_Z/G$. 
Let $\Kum(A_Z)$ be the minimal resolution of $A_Z/G$.
This is a K3 surface over $K$; indeed, it is $\overline{K}$-isomorphic to $\Kum(A)$. Alternatively, we could use the fact that the $A[2]$-action on $A$ gives rise to an action on $\Kum(A)$ and then form the twist $(\Kum(A)\times_K Z)/A[2]$.

Since translations by elements of $A(\ov K)$ act trivially on the cohomology of $A$, we have canonical isomorphisms of $\Ga_K$-modules
$$\H^i_\et(A_{Z,\overline{K}},\Q_\ell)\cong \H^i_\et(A_{\overline{K}},\Q_\ell),\quad i\geq 0.$$ Thus we have a canonical isomorphism
\[ \H^2_\et(\Kum(A_Z)_{\overline{K}},\Q_\ell) \cong \H^2_\et(A_{\overline{K}},\Q_\ell) \oplus \Q_\ell^{Z(\overline{K})}(-1) \]
of $\Ga_K$-representations.

\begin{lem} \label{lem: twisted prelim}
Suppose $\Kum(A_Z)$ has good reduction, and let $\ell$ be an odd prime. Then both $A[2](\overline{K})$ and
$\H^2_\et(A_{\overline{K}},\Q_\ell)$ are unramified as $\Ga_K$-modules, and $[Z]\in   \ker\left( \H^1(K,A[2])\rightarrow \H^1(K_{\nr},A[2])\right)$.
\end{lem}

\begin{proof}
Since $\H^2_\et(\Kum(A_Z)_{\overline{K}},\Q_\ell)$ is an unramified $\Ga_K$-module, we see immediately that both $\H^2_\et(A_{\ov{K}},\Q_\ell)$ and $\Q_\ell^{Z(\ov{K})}(-1)$ are unramified. Thus $Z(\ov K)=Z(K_{\nr})$,
which implies that $A[2](\overline{K})$ is unramified. The same fact also implies that
$Z$ is trivialised over $K_{\nr}$; hence $[Z]$ maps to zero in $\H^1(K_{\nr},A[2])$. 
\end{proof}

It follows that some quadratic twist of $A$ has good reduction (see the appendix).
Since quadratic twists of $A$ do not change the $K$-isomorphism class of $\Kum(A_Z)$, we will therefore assume that $A$ has good reduction, with N\'eron model $\sA/\O_K$. The connected-\'etale sequence for $\sA[2]$ then gives rise to an exact sequence
\[ 0 \lra \sA[2]^\circ_K \lra A[2] \lra \sA[2]^\et_K \lra 0 \]
of group schemes over $K$. Define $Z^\et$ as the pushout of $Z$ along $A[2] \rightarrow \sA[2]^\et_K$; that is, $Z^\et$ is the quotient $Z/\sA[2]^\circ_K$.
We write $\pi\colon Z\to Z^\et$ for the quotient morphism.

\begin{lem} \label{sec}
Suppose that there exists an isomorphism of $K$-schemes $Z\cong \sA[2]^\circ_K \times_K Z^\et$. Then the morphism $\pi\colon Z\rightarrow Z^\et$ has a section.
\end{lem}

\begin{rem} The lemma is \emph{not} immediate since the second projection defined by the isomorphism $Z\cong \sA[2]^\circ_K \times_K Z^\et$ need not coincide with $\pi\colon Z\rightarrow Z^\et$.
\end{rem}

\begin{proof} First of all, let us suppose that $\sA_k$ is ordinary. Write $Z^\et= \Spec(L_1)\sqcup\dots \sqcup \Spec(L_m)$, where each $L_i$ is a finite field extension of $K$. If $L_i=K$ for some $i$, then both $Z$ and $Z^\et$ are trivial torsors, and hence the claim follows from Theorem~\ref{t0}. We may therefore assume either that $m=1$ and $L_1$ is a quartic extension of $K$, or that $m=2$ and $L_1$ and $L_2$ are (not necessarily distinct) quadratic extensions of $K$.

Since $\pi\colon Z\rightarrow Z^\et$ is a $\sA[2]^\circ_K$-torsor, we get a decomposition of $K$-schemes
\[ Z \cong Z_1 \sqcup \dots \sqcup Z_m, \]
where each $Z_i\rightarrow \Spec(L_i)$ is a $\sA[2]^\circ_{K}$-torsor. Now, our given isomorphism $Z\cong \sA[2]^\circ_K \times_K Z^\et$ of $K$-schemes, together with the identity element of $\sA[2]^\circ_K$, gives rise to a $K$-morphism $Z^\et\rightarrow Z$, and hence to $K$-morphisms $\Spec(L_i) \rightarrow Z_{\phi(i)} \rightarrow \Spec(L_{\phi(i)})$ for some function (not necessarily a permutation) $\phi\colon\{1,\ldots,m\}\rightarrow \{1,\ldots,m\}$.

If $\phi$ is a permutation, then the composite map $Z^\et\rightarrow Z \rightarrow Z^\et$ is an isomorphism, and we therefore obtain a section as claimed. But if $\phi$ is not a permutation, then we must have $m=2$, $L_1\cong L_2$ are isomorphic quadratic extensions of $K$, and (after possibly reindexing) $\phi(i)=1$ for $i=1,2$. Write $L=L_1=L_2$. In this case, we have $Z=Z_1\sqcup Z_2$, where $Z_1\rightarrow \Spec(L)$ is a trivial $\sA[2]^\circ_K$-torsor and $Z_2\rightarrow \Spec(L)$ is a possibly non-trivial $\sA[2]^\circ_K$-torsor. But now the fact that we have $K$-isomorphisms
\[ Z = Z_1 \sqcup Z_2 \cong \sA[2]^\circ_K \times_K Z^\et= \sA[2]^\circ_L \sqcup \sA[2]^\circ_L\quad\text{and}\quad
 Z_1\cong \sA[2]^\circ_L\]
implies that there also exists a $K$-isomorphism $Z_2\cong \sA[2]^\circ_L$. This then implies in turn that $Z_2\rightarrow \Spec(L)$ is a trivial $\sA[2]^\circ_K$-torsor. Hence $Z\rightarrow Z^\et$ is a trivial $\sA[2]^\circ_K$-torsor and therefore admits a section as claimed.

The proof in the almost ordinary and supersingular cases is similar, but much easier.
\end{proof}

\bthe \label{t2}
Let $A$ be an abelian surface over $K$ with good, non-supersingular reduction.
Let $Z$ be a $K$-torsor for $A[2]$.
\begin{enumerate}
\item  \label{num: t2 1} The twisted Kummer surface $\Kum(A_Z)$ attains good reduction after a finite extension of $K$.
\item
\label{num: t2 2}
If the \'etale $K$-scheme $Z$ is unramified, then
$\Kum(A_Z)$ attains good reduction after a finite unramified extension of $K$.
\item \label{num: t2 3}
  \begin{enumerate}[label={\rm(\alph*)}, ref={\rm 3\alph*}]
\item\label{num: t2 3 a}  If $A$ has ordinary reduction, then $\Kum(A_Z)$ has good reduction over $K$ if and only if 
the \'etale $K$-scheme $Z^\et$ is unramified and the morphism $\pi\colon Z\rightarrow Z^\et$ has a section.
\item\label{num: t2 3 b} If $A$ has almost ordinary reduction, then $\Kum(A_Z)$ has good reduction over $K$ if and only if the \'etale $K$-scheme $Z^\et$ is unramified, the morphism $\pi\colon Z\rightarrow Z^\et$ has a section, and all points of $\sA[2]^\circ(\overline{K})$ are defined over the splitting field of $Z^\et$.
  \end{enumerate}
\end{enumerate}
\ethe

\begin{rem}\leavevmode
\begin{enumerate}
\item In the almost ordinary case, $Z^\et$ is a $K$-torsor for $\Z/2$, and so its splitting field is either $K$ itself or a quadratic extension of $K$.
\item If $Z$ is a trivial torsor, then the condition in part~\eqref{num: t2 3 a} does not quite reduce to that in Theorem~\ref{t1} since here we only require a scheme-theoretic section, whereas in Theorem~\ref{t1} the section is required to be a group homomorphism. This implicitly proves (in the ordinary case) that the connected-\'etale sequence
\[ 0 \lra \sA^\vee[2](\bar k) \lra A[2](\overline{K})\lra \sA[2](\bar k) \lra 0 \]
has a section as $\Ga_K$-modules if and only if it has a section as $\Ga_K$-sets. The most direct proof that we know of this result goes through Theorem~\ref{t0}.
\item It is not difficult to check that if $\Kum(A_Z)$ has good reduction, then it does so with a scheme model. Indeed, we can easily adapt the explicit constructions of Section~\ref{sec: explicit} to produce scheme-theoretic models of our twisted Kummer surfaces. 
\end{enumerate} 
\end{rem}

\begin{proof}[Proof of Theorem~\ref{t2}]
Parts~\eqref{num: t2 1} and~\eqref{num: t2 2} follow from Theorem~\ref{t1}. Indeed, $\Kum(A)$ and $\Kum(A_Z)$ are isomorphic over a finite extension $L/K$. 
If the $K$-scheme $Z$ is unramified, then $A[2](\overline{K})$ is unramified; hence $L$ can be taken 
to be unramified over $K$.

Now let us turn to part~\eqref{num: t2 3}. We first claim that all hypotheses imply that $Z$ extends to an \'etale $\sA[2]$-torsor $\mathcal{Z}/\O_K$ or, equivalently, that $Z$ is unramified as an \'etale $K$-scheme. Indeed, if $\Kum(A_Z)$ has good reduction, then this follows from Lemma~\ref{lem: twisted prelim}. On the other hand, if $Z^\et$ is unramified and $\pi$ has a section, then the fact that $\pi\colon Z\rightarrow Z^\et$ is a $\sA[2]_K^\circ$-torsor implies that $Z\cong \sA[2]_K^\circ\times_K Z^\et$. In the ordinary case, this implies directly that $Z$ is unramified. In the almost ordinary case, we use the extra assumption that all points of $\sA[2]^\circ(\overline{K})$ are defined over the splitting field of $Z^\et$ to conclude this.

We may therefore assume that we have such a torsor $\mathcal{Z}$. We let $\mathcal{Z}^\et$ denote the pushout of $\mathcal{Z}$ along $\sA[2]\rightarrow \sA[2]^\et$; thus $\mathcal{Z}^\et_K \cong Z^\et$. The special fibres $\mathcal{Z}_k$ and $\mathcal{Z}^\et_k$ are then \'etale torsors for $\sA_k[2]$ and $\sA_k[2]^\et$, respectively. Now, the connected-\'etale sequence for $\sA_k[2]$ splits canonically, and hence we may write $\mathcal{Z}_k=\mathcal{Z}_k^\circ\times_k \mathcal{Z}_k^\et$, where $\mathcal{Z}_k^\circ$ is
an \'etale torsor for $\sA_k[2]^\circ$. Since $\sA_k[2]^\circ$ is trivial as an \'etale sheaf on $\Spec(k)$, any such torsor must be trivial, so we have $\mathcal{Z}_k \cong \sA_k[2]^\circ\times_k \mathcal{Z}^\et_k$. There is a natural action of both $\sA_k[2]$ and $\sA_k[2]^\et$ on $\sA_k$, and we may therefore form the twisted Kummer surfaces $\Kum(\sA_{k,\mathcal{Z}_k})$ and $\Kum(\sA_{k,\mathcal{Z}^\et_k})$, which are isomorphic since $\mathcal{Z}_k \cong \sA_k[2]^\circ\times_k \mathcal{Z}^\et_k$. We claim that $\Kum(\sA_{k,\mathcal{Z}_k})\cong \Kum(\sA_{k,\mathcal{Z}^\et_k})$ is the canonical reduction of $\Kum(A_Z)$. To see this, we can twist the abelian scheme $\sA$ by the $\sA[2]$-torsor $\mathcal{Z}$ to form $\sA_{\mathcal Z}$. Now taking the quotient of $\sA_{\mathcal Z}$ by the natural $G$-action gives rise to a flat $\O_K$-scheme $\sA_{\mathcal Z}/G$, with normal fibres, such that the minimal resolutions of the fibres are $\Kum(A_Z)$ and $\Kum(\sA_{k,\mathcal{Z}_k})$, respectively. Thus we may apply Lemma~\ref{lem: can red alt}. 

As $\Ga_K$-representations, we have
\[\H^2_\et(\Kum(A_Z)_{\overline{K}},\Q_\ell) \cong \H^2_\et(A_{\overline{K}},\Q_\ell) \oplus \Q_\ell^{Z(\overline{K})}(-1) ,  \]
as well as 
\[ \H^2_\et(\Kum(\sA_{k,\mathcal{Z}^\et_k})_{\bar k},\Q_\ell) \cong \H^2_\et(A_{\overline{K}},\Q_\ell) \oplus \Q_\ell^{\sA[2]^\circ(\overline{K}) \times Z^\et(\overline{K})}(-1)\]
if $\sA_k$ is ordinary, and 
\[ \ \H^2_\et(\Kum(\sA_{k,\mathcal{Z}^\et_k})_{\bar k},\Q_\ell) \cong \H^2_\et(A_{\overline{K}},\Q_\ell) \oplus \Q_\ell^{\bigsqcup_{i=1}^8 Z^\et(\overline{K}) }(-1) \]
if $\sA_k$ is almost ordinary. Indeed, in the latter case we know that all sixteen exceptional curves of $\Kum(\sA_{k,\mathcal{Z}^\et_k})$ have to be defined over the (at most) quadratic extension of $k$ trivialising $\mathcal{Z}_k^\et$ but not over any smaller field. We therefore deduce from \cite[Theorem 1.6]{CLL} that $\Kum(A_Z)$ has good reduction over $K$ if and only if there exists an isomorphism of (unramified) $\Ga_K$-sets
\[ Z(\overline{K}) \cong \sA[2]^\circ(\overline{K}) \times Z^\et(\overline{K}) \]
in the ordinary case, and
\[ Z(\overline{K}) \cong \bigsqcup_{i=1}^8 Z^\et(\overline{K}) \]
in the almost ordinary case. 

In the ordinary case, if we have such an isomorphism, then the existence of a section was proved in Lemma~\ref{sec}. Conversely, if $\pi$ has a section, 
then we know that $Z\cong \sA[2]^\circ_K\times_K Z^\et $ since $\pi\colon Z\rightarrow Z^\et$ is a trivial $\sA[2]^\circ_K$-torsor.

In the almost ordinary case, if we have such an isomorphism $Z(\overline{K}) \cong \bigsqcup_{i=1}^8 Z^\et(\overline{K})$,
then \emph{any} set-theoretic section will be $\Ga_K$-equivariant. The fact that all points of $\sA[2]^\circ(\overline{K})$ have to be defined over the splitting field of $Z^\et$ follows from applying Theorem~\ref{t1} over this field.

Conversely, suppose that $\pi$ has a section and all points of $\sA[2]^\circ(\overline{K})$ are defined over a splitting field $L/K$ for $Z^\et$. If $L=K$, then clearly both $Z(\overline{K})$ and $\bigsqcup_{i=1}^8Z^\et(\overline{K})$ are trivial $\Ga_K$-sets. Otherwise, $L/K$ is quadratic, and then both $Z(\overline{K})$ and $\bigsqcup_{i=1}^8Z^\et(\overline{K})$ consist of sixteen points, none of which are fixed by $\Ga_K$, but all of which are fixed by the index $2$ subgroup $\Ga_L\subset \Ga_K$. We can therefore directly construct a $\Ga_K$-equivariant bijection $Z(\overline{K})  \cong \bigsqcup_{i=1}^8Z^\et(\overline{K})$. 
\end{proof}

It is possible for $\Kum(A_Z)$ to have good reduction over $K$ even if $\Kum(A)$ does not. We give examples in both the ordinary and almost ordinary cases. 

\begin{ex} 
To give an example where $A$ has good, ordinary reduction, we take elliptic curves $\mathcal{E}_1$, $\mathcal{E}_2$ over $\Z_2$ with ordinary reduction, such that $\mathcal{E}_1[2](\overline{\Q}_2)$ is a trivial $\Ga_{\Q_2}$-module but $\mathcal{E}_2[2](\overline{\Q}_2)$ is a non-trivial but unramified $\Ga_{\Q_2}$-module. For example, we could take $\mathcal{E}_1$ to be defined by 
\[ y^2+xy=x^3-4x-1  \]
and $\mathcal{E}_2$ to be the curve from Example~\ref{exa: nonsplit}.

Let $K/\Q_2$ be the unramified quadratic extension over which all $\overline{\Q}_2$-points of $\mathcal{E}_2[2]$ are defined (equivalently, over which the connected-\'etale sequence for $\mathcal{E}_2[2]$ splits), and let $\sigma\colon \Z/2 \cong \mathcal{E}_1[2]^\et\rightarrow \mathcal{E}_1[2]$ be a splitting of the connected-\'etale sequence of $\mathcal{E}_1$. We take $\sA=\mathcal{E}_1\times_{\Z_2}\mathcal{E}_2$ and $A=\sA_{\Q_2}$; thus $\sigma$ induces a map $\Z/2\rightarrow A[2]$. The quadratic extension $K/\Q_2$ gives rise to a class $[K]\in \H^1(\Ga_{\Q_2},\Z/2)$, and we let $Z/\Q_2$ be an $A[2]$-torsor whose cohomology class is equal to the image of $[K]$ in $\H^1(\Ga_{\Q_2},A[2])$. 

We then see that $Z$ and $\sA[2]^\circ_{\Q_2} \times_{\Q_2} Z^\et$ are isomorphic as $\Q_2$-schemes since neither admits a $\Q_2$-point but both split into sixteen disjoint rational points over the quadratic extension $K/\Q_2$. Hence the twisted Kummer surface $\Kum(A_Z)$ has good reduction by Theorem~\ref{t2} and Lemma~\ref{sec}. On the other hand, since the connected \'etale sequence for $\mathcal{E}_2$ is not split, neither is that for $\sA$, and hence the Kummer surface $\Kum(A)$ cannot have good reduction by Theorem~\ref{t1}.
\end{ex}

\begin{ex} For an example in the almost ordinary case, we choose a $2$-adic field $K$ and elliptic curves $\mathcal{E}_1,\mathcal{E}_2/\O_K$ such that $\mathcal{E}_1$ has ordinary reduction, $\mathcal{E}_2$ has supersingular reduction, $\mathcal{E}_1$ has full $2$-torsion defined over $K$, but $\mathcal{E}_2$ has one $K$-point of order $2$, with the remaining two defined over an unramified quadratic extension $L/K$. 

To see that such examples exist, we take $K=\Q_2(\sqrt[3]{2})$, take $\mathcal{E}_1/\O_K$ defined by
\[  y^2+xy=x^3-\frac{81}{389017}, \]
and take $\mathcal{E}_2/\O_K$ defined by
\[ y^2+y=x^3. \]
The $2$-torsion of $\mathcal{E}_2$ is defined over the unramified quadratic extension $L=K(\omega)$, where $\omega$ is a primitive third root of unity.

Now we take $\sA=\mathcal{E}_1\times_{\O_K}\mathcal{E}_2$. This has almost ordinary reduction, but $A:=\sA_K$ does not have full $2$-torsion defined over $K$. Hence $\Kum(A)$ does not have good reduction over $K$ by Theorem~\ref{t0}. On the other hand, if we let $\mathcal{Z}^\et/\O_K$ be the unique non-trivial $\sA[2]^\et$-torsor trivialised by the unramified extension $L/K$ and choose a splitting $\sigma\colon \sA[2]^\et\rightarrow \sA[2]$ of the connected \'etale sequence (which exists because $\mathcal{E}_1$ has full $2$-torsion defined over $K$), then we may take $Z$ to be the generic fibre of 
the $\sA[2]$-torsor over $\O_K$ which is the pushout of $\mathcal{Z}^\et$ under the map
$\sigma\colon \sA[2]^\et\rightarrow \sA[2]$.
It is straightforward to check that $Z$ and $\sA$ satisfy the hypotheses of Theorem~\ref{t2}\eqref{num: t1 3 b};
thus $\Kum(A_Z)$ has good reduction.
\end{ex}

\renewcommand\thesection{\Alph{section}}
\setcounter{section}{0}

\section*{Appendix. Good reduction and quadratic twists} \label{app: twists}
\addcontentsline{toc}{section}{Appendix. Good reduction and quadratic twists}
\refstepcounter{section}

In this appendix we prove the following result, which is surely well known to the experts. For this appendix only, we drop the assumption that $\Char(k)=2$; thus $\O_K$ will be a complete DVR with fraction field $K$ of characteristic $0$ and perfect residue field $k$ of arbitrary characteristic.

\bthe \label{theo: twist 1}
Let $A$ be an abelian surface over $K$, and suppose that $\Kum(A)$
has good reduction over $K$. Then there exists a character $\chi\colon\Ga_K\to\Z/2$ such that the quadratic twist $A^\chi$ has good reduction over $K$.
\ethe

Combining the smooth and proper base change theorem with the classical N\'eron--Ogg--Shafarevich criterion, and using the inclusion $\H^2_\et(A_{\overline{K}},\Q_\ell)\subset \H^2_\et(\Kum(A)_{\overline{K}},\Q_\ell)$, we can restate Theorem~\ref{theo: twist 1} as follows.

\bthe \label{theo: twist 2} Suppose that $\ell\neq \Char(k)$ and that the inertia subgroup ${\rm I}_K\subset \Ga_K$ acts trivially on $\H^2_\et(A_{\overline{K}},\Q_\ell)$. Then there exists 
a character $\chi\colon\Ga_K\to\Z/2$ such that ${\rm I}_K$ acts trivially on $\H^1_\et(A^\chi_{\overline{K}},\Q_\ell)$.
\ethe
\begin{proof}
Let $\rho\colon \Ga_K\to \GL(\H^1_\et(A_{\overline{K}},\Q_\ell))$ be the representation
of $\Ga_K$ in $\H^1_\et(A_{\overline{K}},\Q_\ell)$.
The triviality of the action of ${\rm I}_K$ on $\H^2_\et(A_{\overline{K}},\Q_\ell)$ implies 
that $\rho({\rm I}_K)\subset\{\pm 1\}$; see \cite[Lemma 4.4(ii)]{Ove21}.
Since $\{\pm 1\}$ is central in $\GL(\H^1_\et(A_{\overline{K}},\Q_\ell))$,
the restriction of $\rho$ to ${\rm I}_K$ is invariant under the action of $\Ga_K$ on ${\rm I}_K$
by conjugation. This implies that the restriction of $\rho$ to ${\rm I}_K$ is a $\Ga_k$-invariant element of
$\H^1({\rm I}_K,\Z/2)=\Hom({\rm I}_K,\Z/2)$.

We next claim that the restriction map 
\[ \H^1(\Ga_K,\Z/2)\lra \H^1({\rm I}_K,\Z/2)^{\Ga_k} \]
is surjective. By the Hochschild--Serre spectral sequence, this is equivalent to injectivity of the inflation map
\[ \H^2(\Ga_k,\Z/2)\lra \H^2(\Ga_K,\Z/2). \]
If $\Char(k)\neq 2$, this is well known; see, \textit{e.g.}, \cite[Equation~(1.15), p.~30]{CTS21}.
If $\Char(k)=2$, then the Artin--Schreier exact sequence implies that $\H^2(\Ga_k,\Z/2)=0$; see,
\textit{e.g.}, \cite[\href{https://stacks.math.columbia.edu/tag/0A3K}{Lemma 0A3K}]{stacks-project}. Hence the map is trivially injective. 
 
To complete the proof, we now choose $\chi\in\H^1(\Ga_K,\Z/2)=\Hom(\Ga_K,\Z/2)$ 
to be a lifting of the restriction of $\rho$ to ${\rm I}_K$. Then the representation of $\Ga_K$
in $\H^1_\et(A^\chi_{\overline{K}},\Q_\ell)$ restricts to the trivial representation of ${\rm I}_K$.
\end{proof}

\newcommand{\etalchar}[1]{$^{#1}$}
\providecommand{\bysame}{\leavevmode\hbox to3em{\hrulefill}\thinspace}

\end{document}